\documentclass[12pt]{article}
\usepackage{amssymb,amsmath,amsfonts}
\usepackage{color}
\usepackage[margin=0.8in]{geometry}
\usepackage[T1]{fontenc}

\newcommand{\dd}{\,\mathrm{d}}
\newcommand{\bg}{\mathbf{g}}
\newcommand{\bh}{\mathbf{h}}
\newcommand{\bA}{\mathbf{A}}
\newcommand{\bs}{\mathbf{S}}
\newcommand{\bu}{\mathbf{u}}
\newcommand{\bv}{\mathbf{v}}
\newcommand{\bw}{\mathbf{w}}
\newcommand{\bx}{\mathbf{x}}
\newcommand{\bF}{\mathbf{f}}
\newcommand{\bG}{\mathbf{b}}
\newcommand{\bT}{\mathbf{T}}
\newcommand{\bI}{\mathbf{I}}
\newcommand{\bsig}{\boldsymbol{\sigma}}
\newcommand{\beps}{\boldsymbol{\varepsilon}}
\newcommand{\bkap}{\boldsymbol{\kappa}}
\newcommand{\bphi}{\boldsymbol{\varphi}}
\newcommand{\CC}{\mathrm{C}}
\newcommand{\LL}{\mathrm{L}}
\newcommand{\HH}{\mathrm{H}}
\newcommand{\WW}{\mathrm{W}}
\DeclareMathOperator{\tr}{tr}
\DeclareMathOperator{\Div}{Div}

\newtheorem{definition}{Definition}
\newtheorem{theorem}{Theorem}
\newtheorem{lemma}{Lemma}
\newtheorem{remark}{Remark}

\numberwithin{equation}{section}
\numberwithin{theorem}{section}
\numberwithin{lemma}{section}
\numberwithin{definition}{section}
\numberwithin{remark}{section}

\begin{document}

\title{\textbf{Well-posedness of the fractional Zener wave equation for heterogenous viscoelastic materials}}
\author{\sc{Ljubica Oparnica and Endre S\"uli}\\~\\}

\date{~}
\maketitle

\begin{abstract}
{\color{black} The Zener model for viscoelastic solids replaces Hooke's law} $\bsig = 2\mu \beps(\bu) + \lambda \tr(\beps(\bu)) \bI$, relating the stress tensor $\bsig$ to the strain tensor $\beps(\bu)$, where $\bu$ is the displacement vector, $\mu>0$ is the shear modulus, and $\lambda\geq 0$ is the first Lam\'{e} coefficient,  with the constitutive law
$
(1 + \tau D_t) \bsig = (1 + \rho D_t)[2\mu \beps(\bu) + \lambda \tr(\beps(\bu)) \bI]$,
where $\tau>0$ is the characteristic relaxation time and $\rho \geq \tau$ is the characteristic retardation time. It is the simplest model that predicts creep/recovery and stress relaxation phenomena. We explore the well-posedness of the {\color{black} fractional version of the} model, where the first-order time-derivative $D_t$ in the constitutive law is replaced by the  {\color{black} Caputo} time-derivative $D_t^\alpha$, with $\alpha \in (0,1)$, $\mu, \lambda$ belong to $\LL^\infty(\Omega)$, $\mu$ is bounded below by a positive constant and $\lambda$ is nonnegative. We show that, when coupled with the equation of motion $\varrho \ddot{\bu} = \Div \bsig + \bF$, considered in a bounded open Lipschitz domain $\Omega$ in $\mathbb{R}^3$ and over a time interval $(0,T]$, where $\varrho \in \LL^\infty(\Omega)$ is the density of the material, assumed to be bounded below by a positive constant, and $\bF$ is a specified load vector, the resulting model is well-posed in the sense that the associated initial-boundary-value problem, {\color{black} with initial conditions $\bu(0,\bx) = \bg(\bx)$, $\dot{\bu}(0,\bx) = \bh(\bx)$, $\bsig(0,\bx) = \bs(\bx)$, for $\bx \in \Omega$, and a homogeneous Dirichlet boundary condition, possesses a unique weak solution for any choice of $\bg \in [\HH^1_0(\Omega)]^3$, $\bh \in [\LL^2(\Omega)]^3$, and $\bs = \bs^{\rm T} \in [\LL^2(\Omega)]^{3 \times 3}$, and any load vector $\bF \in \LL^2(0,T;[\LL^2(\Omega)]^3)$, and that this unique weak solution depends continuously on the initial data and the load vector.}
\end{abstract}

\section{Statement of the model}
Suppose that $\Omega \subset \mathbb{R}^3$ is a bounded, open, simply-connected Lipschitz domain,
with boundary $\partial \Omega$, occupied by a viscoelastic material, and let $T>0$. Consider the equation of motion
\begin{align}\label{eq:1}
\varrho \ddot{\bu} = \Div{\bsig} + \bF\qquad \mbox{in $(0,T] \times \Omega$},
\end{align}
with $\varrho>0$ signifying the density of the material,
$\bu$ the displacement vector, $\bsig$ the stress tensor, and $\bF$ the load vector,
with the material being considered subject to the initial conditions
\begin{align}\label{eq:2}
\bu(0,\bx) = \bg(\bx),\qquad \dot{\bu}(0,\bx) = \bh(\bx),\qquad \bsig(0,\bx) = \bs(\bx),\qquad \mbox{for $\bx \in \Omega$},
\end{align}
and a suitable boundary condition, which
for the sake of simplicity of the exposition we shall assume to be the homogeneous Dirichlet boundary condition
\begin{align}\label{eq:3}
\bu(t,\bx) = \mathbf{0} \qquad \mbox{for all $(t,\bx) \in (0,T] \times \partial \Omega$}.
\end{align}
The discussion below trivially extends to the case of a {\color{black} mixed homogeneous Dirichlet/             nonhomoge\-neous} Neumann boundary condition provided that the Dirichlet part of $\partial\Omega$ has positive two-dimensional surface measure {\color{black} (cf. the concluding remarks at the end of the paper for further comments in this direction)}. In the case of a classical linear (Hookean) elastic body the stress tensor $\bsig$ is related to the strain tensor (symmetric displacement gradient)
\[ \beps(\bu):= \frac{1}{2}(\nabla \bu + (\nabla \bu)^{\rm T}) \]
through \textit{Hooke's law}
$\bsig = 2\mu \beps(\bu) + \lambda \tr(\beps(\bu)) \bI$,
where $\mu>0$ is the \textit{shear modulus} and $\lambda\geq 0$ is the \textit{first Lam\'{e} coefficient}. In this case the initial value $\bs=\bsig|_{t=0}$ of $\bsig$ is automatically equal to $2\mu \beps(\bg) + \lambda \tr(\beps(\bg)) \bI$, by Hooke's law, and need not (or, more precisely, should not) be specified independently, as otherwise the resulting initial-boundary-value problem will be over-determined and will have no solution in general. However for Zener's model under consideration here the situation is different: the constitutive law relating the stress tensor $\bsig$ to the strain tensor $\beps(\bu)$ involves the time-derivative of order $\alpha \in (0,1]$ of $\bsig$:
\begin{align*}
(1 + \tau^\alpha D^\alpha_t) \bsig = (1 + \rho^\alpha D^\alpha_t)[2\mu \beps(\bu) + \lambda \tr(\beps(\bu)) \bI],
\end{align*}
with $\tau>0$ signifying the \textit{characteristic relaxation time} and $\rho \geq \tau$ is the \textit{characteristic
retardation time}, --- which then necessitates the specification of an initial datum $\bs$ for $\bsig$.

In the case of $\alpha=1$ the model was proposed by Zener \cite{Zener}
(with $\lambda = 0$). The fractional version of Zener's model was
introduced (in one space dimension and, again, with $\lambda=0$) by
Caputo and Mainardi (cf. \cite{CM}, and eq. (13) in \cite{FD}); in the
context of the present paper a natural generalization of the model
from \cite{CM} to the case of three space-dimensions would be
\[ (1 + \tau^\alpha D^\alpha_t) \bsig = E (1  +\, \rho^\alpha
D^\alpha_t) \beps(\bu), \quad \mbox{with $\bsig(0,\cdot) = E
	\left(\frac{\rho}{\tau}\right)^\alpha \beps(\bu(0,\cdot))$},\]
where, following Bagley and Torvik \cite{BT},  $E>0$ is referred to as
the \textit{rubbery modulus}, $E(\rho/\tau)^\alpha$ is called the
\textit{glassy modulus}, and $\alpha \in (0,1)$ is the
\textit{fractional order of evolution}.
As has been noted by Freed and Diethelm \cite{FD}, this model
allows for a finite discontinuity in the stress-strain response at
time zero (cf. Remark \ref{re:1} below for further comments on this observation in the
context of our well-posedness analysis).
Bagley and Torvik \cite{BT} have demonstrated that the fractional orders of
evolution in stress and strain must be the same, as originally proposed in the work of
Caputo and Mainardi \cite{CM}, in order that a material model of fractional order
comply with the second law of thermodynamics; Bagley and Calico
\cite{BC} have also shown that the differential orders need to be the same for
the stress and the strain in order to ensure that sound waves in the material propagate at finite
speed. For further motivation from the point of view of continuum thermodynamics
for considering fractional-order constitutive laws of this kind we refer
to \cite{A}, \cite{BC}, \cite{BT}, and \cite{P}, for example.

As the actual value of the characteristic retardation time $\rho ~(\geq \tau > 0)$ is of no relevance in the discussion that follows, for the sake of simplicity of the exposition we have fixed $\rho=1$, resulting in the constitutive law
\begin{align}\label{eq:4}
(1 + \tau^\alpha D^\alpha_t) \bsig = (1 + D^\alpha_t)[2\mu \beps(\bu) + \lambda \tr(\beps(\bu)) \bI], \qquad \tau \in (0,1], \quad \alpha \in (0,1).
\end{align}
As will be seen in what follows, the relation $(1=) \rho \geq \tau>0$ is crucial for ensuring the well-posedness of the
resulting model, in agreement with the discussion in \cite{BT} (particularly eqs. (14) and (22)--(25) therein with $\alpha=\beta$) concerning the relevant thermodynamical conditions to ensure nonnegativity of the internal work and guarantee a nonnegative rate of energy dissipation.
The constitutive law \eqref{eq:4} generalizes the one proposed by Caputo and Mainardi in
\cite{CM} in that we admit $\lambda \geq 0$, motivated by the fact that formally setting $\alpha=0$ in \eqref{eq:4} reduces it to Hooke's constitutive law. As a matter of fact, we shall assume, more generally, that
\begin{alignat}{2}\label{coeff-ass}
\begin{aligned}
\varrho \in \LL^\infty(\Omega), \qquad \mbox{and there exists a positive constant $\varrho_0$ such that $\varrho(\bx) \geq \varrho_0$ a.e. in $\Omega$},\\
\mu \in \LL^\infty(\Omega), \qquad \mbox{and there exists a positive constant $\mu_0$ such that $\mu(\bx) \geq \mu_0$ a.e. in $\Omega$},\\
\lambda \in \LL^\infty(\Omega), \qquad \mbox{and $\lambda(\bx) \geq 0$ a.e. in $\Omega$},
\end{aligned}
\end{alignat}
so as to admit spatially heterogeneous viscoelastic materials.
With straightforward modifications all of our results extend to the
case of Hooke's model corresponding to $\alpha=0$ and the classical Zener model corresponding to $\alpha=1$; we shall therefore confine ourselves to the, technically more involved, fractional-order setting, when $\alpha \in (0,1)$.

Zener's constitutive law aims to overcome some of the shortcomings of the Maxwell and Kelvin--Voigt models: the Maxwell model does not describe creep or recovery, and the Kelvin--Voigt model does not describe stress relaxation. Zener's constitutive law is the simplest model that predicts both phenomena. Our aim here is to explore the well-posedness of the model, focusing in particular on its refinement, where the first time-derivative $D_t$ featuring in the constitutive law is replaced by a fractional-order time-derivative $D_t^\alpha$, with $\alpha \in (0,1)$. We emphasize that the equation of motion \eqref{eq:1}, expressing balance of the linear momentum in terms of the Cauchy stress, remains unchanged: it is
only the constitutive law relating the stress tensor to the strain tensor, which encodes the specific properties of the material, that is altered here by admitting the fractional range $\alpha \in (0,1)$.

The fractional derivative $D^\alpha_t$ of order $\alpha \in (0,1)$ appearing in \eqref{eq:4}
is in the sense of Caputo. It is understood to be acting on $3$-component vector-functions and $3 \times 3$-matrix-valued functions componentwise. In particular, for a scalar-valued function $f \in \mbox{AC}([0,T])$,
\[ (D^\alpha_t f)(t):= \frac{1}{\Gamma(1-\alpha)} \int_0^t \frac{\dot{f}(s)}{(t-s)^\alpha} \dd s,\qquad t \in (0,T].\]

{\color{black}
The partial differential equation \eqref{eq:1} coupled with the constitutive law \eqref{eq:4} is referred to as the
\textit{fractional Zener wave equation}. Wave propagation in viscoelastic media governed by the fractional Zener constitutive law in one space dimension was first considered by Caputo and Mainardi \cite{CM}. The existence and uniqueness of  the fundamental solution of a generalized Cauchy problem for the fractional Zener wave equation were proved in \cite{KOZ2010}, and an explicit expression for the solution was also given (cf. Theorem 4.2 in \cite{KOZ2010}).
The existence and uniqueness of solutions for a generalization of the fractional Zener wave equation
proposed by Enelund and Josefson \cite{EJ},
in the case of mixed homogenous Dirichlet/nonhomogeneous Neumann boundary conditions on bounded polytopal domains in two
and three space dimensions, were proved by Saedpanah in \cite{FS2014}; and, under suitable restrictions on the domain $\Omega$ and the data, weak solutions of the model were shown in \cite{FS2014} to possess additional regularity. In an earlier work,
Larsson and Saedpanah \cite{LS} showed the well-posedness of the homogeneous Dirichlet problem for this model using techniques from linear semigroup theory. The weak formulation of the evolution equation \eqref{eq:8} that we study here differs from the one considered in \cite{FS2014}; indeed, equation (2.7)$_1$ in \cite{FS2014} was arrived at by using Laplace transform techniques on the constitutive law to obtain an explicit expression for the stress tensor in terms of the strain tensor, which was then substituted into the equation of motion to eliminate the stress tensor; whereas, as we shall explain below, we Laplace transform the equation of motion as well as the constitutive law and we then eliminate the Laplace transform of the stress tensor from the transformed equation of motion.
Furthermore, in both \cite{LS} and \cite{FS2014} the fractional derivative featuring in the constitutive law was the left Riemann--Liouville derivative rather than the Caputo derivative considered here, and the initial response for the
stress tensor was assumed to follow Hooke's law.

The aim of the present work is to explore the question of existence and uniqueness of weak solutions to the initial-boundary-value problem \eqref{eq:1}--\eqref{eq:4} without the additional assumption that the initial response for the stress follows Hooke's law. In the absence of this extra assumption on the initial stress the analysis of the model is considerably more complicated; nevertheless, we are able to show (cf. Theorem \ref{th:1} below) that the model \eqref{eq:1}--\eqref{eq:4} admits a unique weak solution for any $\bF \in \LL^2(0,T;[\LL^2(\Omega)]^3)$, and
arbitrary initial data $\bg \in [\HH^1_0(\Omega)]^3$, $\bh \in [\LL^2(\Omega)]^3$, and $\bs = \bs^{\rm T} \in [\LL^2(\Omega)]^{3 \times 3}$, without any additional restrictions on the choice of $\bs$.

To this end, our first objective is to transform the fractional Zener model \eqref{eq:1}--\eqref{eq:4} to a form in which it is amenable to mathematical analysis. We shall therefore Laplace-transform the equation of motion \eqref{eq:1} (where it will be understood that the source term $\bF$ is extended by $\mathbf{0}$ from $(0,T] \times \Omega$ to $(0,\infty)\times \Omega$), as well as the constitutive law \eqref{eq:4} with respect to the temporal variable $t$ (again with the understanding that, for the moment, $t \in (0,\infty)$ rather than $t\in (0,T]$ with $T<\infty$). This will enable
us to eliminate the stress tensor $\bsig$ from the equation of motion in terms of the strain tensor $\beps(\bu)$, resulting in a second-order nonlocal evolution equation (cf. \eqref{eq:8} below), which will then be the focus of our subsequent analysis. We shall concentrate on the proof of existence and uniqueness of weak solutions, and the continuous dependence of weak solutions on the data. Specifically, we shall show that the constitutive law \eqref{eq:4}, when coupled with \eqref{eq:1}--\eqref{eq:3}, gives rise to a well-posed
mathematical model: by using a compactness argument we shall prove the existence of a weak solution to the model and will prove that weak solutions thus constructed satisfy an energy inequality, which bounds appropriate norms of the solution in terms of norms of the initial data and the source term; we shall also show that weak solutions are unique.
}

\section{Zener's model as a fractional evolution equation}
The aim of this section is to merge the equation of motion \eqref{eq:1} and the constitutive law \eqref{eq:4} into a single evolution equation, which we shall then subject to mathematical analysis. We proceed by eliminating the stress tensor $\bsig$ from
\eqref{eq:1} by Laplace transforming both \eqref{eq:1} and the constitutive law \eqref{eq:4}.

The Laplace transform with respect to the variable $t$ of a function $f$ defined on $(0,\infty)$ such that
$\int_0^\infty  |f(t)|\, \mathrm{e}^{-at} \dd t < \infty$ for some $a \in \mathbb{R}$, is defined by
\[ \mathcal{L}(f)(p) = \tilde{f}(p) := \int_0^\infty f(t)\, \mathrm{e}^{-pt} \dd t,\qquad \mbox{for $p \in \mathbb{C}$ with $\mathrm{Re}\,p \geq a$}.\]

Then, for any $f \in \CC([0,\infty))\cap \CC^1((0,\infty))$ such that $\int_0^\infty (|\dot{f}(t)| + |f(t)|)\, \mathrm{e}^{-at} \dd t< \infty$ for some $a \in \mathbb{R}$, straightforward calculations yield that
\begin{align*}
\mathcal{L}(\dot{f})(p) &= p \tilde{f}(p) - f(0),\qquad \mathrm{Re}\,p \geq a,
\end{align*}
where the symbol $\cdot$ over a $t$-dependent function denotes its derivative with respect to $t$, and, similarly, $\cdot\cdot$ over a $t$-dependent function denotes its second derivative with respect to $t$. As
\[\mathcal{L}((\cdot)^{-\alpha})(p) = \Gamma(1-\alpha) p^{\alpha - 1},\qquad \!\mathrm{Re}\,p > 0,\quad \alpha \in (0,1),\]
by noting that
\[ D^\alpha_t f = \frac{1}{\Gamma(1-\alpha)} \left[ \dot f \ast_t (\cdot)^{-\alpha}\right], \]
where the convolution $\ast_t$ is defined by $(f\ast_t g)(t):= \int_0^t f(s) g(t-s) \dd s$, we have that
\begin{align*}
\mathcal{L}(D^\alpha_t f)(p) &= \frac{1}{\Gamma(1-\alpha)} \mathcal{L}\left[ \dot{f} \ast_t (\cdot)^{-\alpha}\right](p)
= \frac{1}{\Gamma(1-\alpha)} \mathcal{L}(\dot{f})(p)\, \mathcal{L}({(\cdot)}^{-\alpha})(p)\\
&= p^\alpha \tilde{f}(p) - p^{\alpha-1} f(0),\qquad \mathrm{Re}\,p \geq a, \quad \alpha \in (0,1).
\end{align*}

Consider the Mittag--Leffler function
\[ E_{\alpha,\beta}(z):= \sum_{k=0}^\infty \frac{z^k}{\Gamma(\alpha k + \beta)}, \qquad z \in \mathbb{C}, \quad \alpha>0, \quad \beta >0. \]
Letting
\[ e_\alpha(t,\gamma):= E_{\alpha,1}(-\gamma t^\alpha),\qquad t \in [0,\infty),\quad \gamma>0,\]
one has that
\begin{align}\label{eq:5} \mathcal{L}(e_\alpha(\cdot,\gamma))(p) = \frac{p^{\alpha-1}}{p^\alpha + \gamma} \qquad \mbox{for $\mathrm{Re}\, p > \gamma^{\frac{1}{\alpha}}$.}
\end{align}

Henceforth, for the sake of simplicity, we shall write $e_{\alpha,\gamma}(t)$ instead of $e_\alpha(t,\gamma)$, and restrict
ourself to the range $\alpha \in (0,1)$ of relevance to us in the present context. As $e_{\alpha,\gamma}(0)=1$, it follows that
\[ \mathcal{L}(\dot{e}_{\alpha,\gamma})(p) = p\, \tilde{e}_{\alpha,\gamma}(p) - 1 =  \frac{p^{\alpha}}{p^\alpha + \gamma} - \mathcal{L}(\delta),\qquad \mathrm{Re}\, p > \gamma^{\frac{1}{\alpha}},\]
where $\delta$ is the Dirac distribution concentrated at $t=0$. Thus, now with the Laplace transform acting in the sense of tempered distributions\footnote{For a tempered distribution $f \in \mathcal{S}'$, with $\mbox{supp} (f) \subset [0,\infty)$, we define
$\mathcal{L}(f)(p) = \tilde{f}(p) := \langle f  ,  \eta\, \mathrm{e}^{-p\cdot} \rangle$, for $\mbox{Re } p > 0$, where $\eta \in C^\infty(\mathbb{R})$ is such that $\eta(t) \equiv 0$ for $t \leq -2$ and $\eta(t)\equiv 1$ for $t \geq -1$.}
\[ \mathcal{L}(\dot{e}_{\alpha,\gamma} + \delta)(p) =  \frac{p^{\alpha}}{p^\alpha + \gamma}, \qquad \mathrm{Re}\, p > \gamma^{\frac{1}{\alpha}}.\]
As a consequence of this identity we have that
\begin{align}\label{eq:6}
\begin{aligned}
\mathcal{L}^{-1}\left(\frac{1 + \tau^\alpha p^\alpha}{1 + p^\alpha}\right) &= \mathcal{L}^{-1}\left(1 + (\tau^\alpha - 1) \frac{p^\alpha}{1+  p^\alpha}\right) = \delta + (\tau^\alpha-1)\mathcal{L}^{-1}\left(\frac{p^\alpha}{p^\alpha + 1}\right)\\
& = \delta + (\tau^\alpha-1) (\dot{e}_{\alpha,1} + \delta).
\end{aligned}
\end{align}

Following these preparatory considerations, we Laplace-transform the constitutive law \eqref{eq:4}, which yields
\begin{align*}
\tilde{\bsig} + \tau^\alpha \mathcal{L}(D^\alpha_t \bsig) = 2\mu \beps(\tilde\bu) + \lambda \tr(\beps(\tilde\bu)) \bI\, +\, 2\mu  \mathcal{L}(D^\alpha_t \beps(\bu)) + \lambda \mathcal{L}(D^\alpha_t (\tr(\beps(\bu)) \bI),\qquad \tau \in (0,1].
\end{align*}
Hence,
\begin{align*}
&\tilde{\bsig}(p) + \tau^\alpha (p^\alpha \tilde\bsig(p) - p^{\alpha-1}\bs) \\\quad &= 2\mu \beps(\tilde\bu) + \lambda \tr(\beps(\tilde\bu)) \bI\, +\, 2\mu (p^\alpha \beps(\tilde{\bu}) - p^{\alpha-1}\beps(\bg)) + \lambda (p^\alpha \tr(\beps(\tilde\bu)) \bI - p^{\alpha-1} \tr(\beps(\bg)) \bI).
\end{align*}
Equivalently,
\begin{align*}
&(1  + \tau^\alpha p^\alpha) \tilde\bsig(p) = (1+ p^\alpha)(2\mu \beps(\tilde\bu) + \lambda \tr(\beps(\tilde\bu)) \bI)\, +\, p^{\alpha-1} (\tau^\alpha\bs-
2\mu \beps(\bg) - \lambda \tr(\beps(\bg)) \bI),
\end{align*}
and therefore
\begin{align}\label{eq:sig-laplace}
&\tilde\bsig(p) = \frac{1+ p^\alpha}{1  + \tau^\alpha p^\alpha}(2\mu \beps(\tilde\bu) + \lambda \tr(\beps(\tilde\bu)) \bI)\, +\, \frac{p^{\alpha-1}}{1  + \tau^\alpha p^\alpha} (\tau^\alpha\bs-
2\mu \beps(\bg) - \lambda \tr(\beps(\bg)) \bI).
\end{align}

Consequently, and by Laplace-transforming the equation of motion \eqref{eq:1}, we deduce that
\[ \varrho \mathcal{L}(\ddot{\bu}) = \frac{1+ p^\alpha}{1  + \tau^\alpha p^\alpha}\Div(2\mu \beps(\tilde\bu) + \lambda \tr(\beps(\tilde\bu)) \bI)\, +\, \frac{p^{\alpha-1}}{1  + \tau^\alpha p^\alpha} \Div (\tau^\alpha\bs-
2\mu \beps(\bg) - \lambda \tr(\beps(\bg)) \bI) + \tilde{\bF},\]
and, upon multiplying this equality by $\frac{1+ \tau^\alpha p^\alpha}{1  + p^\alpha}$, we have that
\begin{align*}
 &\varrho \frac{1+ \tau^\alpha p^\alpha}{1  + p^\alpha}\mathcal{L}(\ddot{\bu}) = \Div(2\mu \beps(\tilde\bu) + \lambda \tr(\beps(\tilde\bu)) \bI)\,\\
 &\qquad  +\, \frac{p^{\alpha-1}}{1  + p^\alpha} \Div (\tau^\alpha\bs-
2\mu \beps(\bg) - \lambda \tr(\beps(\bg)) \bI) + \frac{1+ \tau^\alpha p^\alpha}{1  + p^\alpha} \tilde{\bF}.
\end{align*}

Hence, by inverse-Laplace-transforming this equality and applying the convolution theorem for the Laplace transform, we obtain
\begin{align*}
 &\varrho \mathcal{L}^{-1}\left(\frac{1+ \tau^\alpha p^\alpha}{1  + p^\alpha}\right) \ast_t \ddot{\bu} = \Div(2\mu \beps(\bu) + \lambda \tr(\beps(\bu)) \bI)\, \\
 &\qquad +\, \mathcal{L}^{-1}\left(\frac{p^{\alpha-1}}{1  + p^\alpha}\right) \Div (\tau^\alpha\bs-
2\mu \beps(\bg) - \lambda \tr(\beps(\bg)) \bI) +  \mathcal{L}^{-1}\left(\frac{1+ \tau^\alpha p^\alpha}{1  + p^\alpha}\right) \ast_t \bF.
\end{align*}
Using \eqref{eq:6} and \eqref{eq:5} we then deduce that
\begin{align*}
\varrho (\delta + (\tau^\alpha-1) \left(\dot{e}_{\alpha,1} + \delta\right))\ast_t \ddot{\bu} &= \Div(2\mu \beps(\bu) + \lambda \tr(\beps(\bu)) \bI)\, \\
&\quad +\, e_{\alpha,1}\, \Div (\tau^\alpha\bs-
2\mu \beps(\bg) - \lambda \tr(\beps(\bg)) \bI)\\
&\quad + (\delta + (\tau^\alpha-1) \left(\dot{e}_{\alpha,1} + \delta\right))\ast_t \bF,
\end{align*}
and therefore
\begin{align*}
\varrho \tau^\alpha \ddot \bu + \varrho \,(\tau^\alpha-1) \dot{e}_{\alpha,1} \ast_t \ddot{\bu} &= \Div(2\mu \beps(\bu) + \lambda \tr(\beps(\bu)) \bI)\\
&\quad +\, e_{\alpha,1}\, \Div (\tau^\alpha\bs-
2\mu \beps(\bg) - \lambda \tr(\beps(\bg)) \bI)\\
&\quad +  \tau^\alpha \bF + (\tau^\alpha-1) \dot{e}_{\alpha,1} \ast_t \bF.
\end{align*}
We now focus on the second term on the left-hand side of this equality. By noting that
\[ (f \ast_t \dot{g})(t) = \frac{\dd}{\dd t}(f\ast_t g)(t)- f(t)g(0)\]
we deduce (by suppressing the $\bx$-dependence of $\bu$ for the sake of notational simplicity) that
\[ (\dot{e}_{\alpha,1} \ast_t \ddot{\bu})(t) = \frac{\partial}{\partial t}(\dot{e}_{\alpha,1}\ast_t \dot{\bu})(t)- \dot{e}_{\alpha,1}(t) \dot{\bu}(0) = \frac{\partial}{\partial t}(\dot{e}_{\alpha,1}\ast_t \dot{\bu})(t)- \dot{e}_{\alpha,1}(t) \bh.\]
Consequently,
\begin{align*}
&\varrho \tau^\alpha \ddot \bu + \varrho (\tau^\alpha-1)\left[\frac{\partial}{\partial t}(\dot{e}_{\alpha,1}\ast_t \dot{\bu})- \dot{e}_{\alpha,1} \bh \right]  \\
&\qquad = \Div(2\mu \beps(\bu) + \lambda \tr(\beps(\bu)) \bI)\,
+\, e_{\alpha,1}\, \Div (\tau^\alpha\bs-
2\mu \beps(\bg) - \lambda \tr(\beps(\bg)) \bI)\\
&\qquad \quad +  \tau^\alpha \bF + (\tau^\alpha-1) \dot{e}_{\alpha,1} \ast_t \bF,
\end{align*}
which upon rearrangement yields
\begin{align}\label{eq:7}
\begin{aligned}
&\tau^\alpha\, \varrho\ddot \bu + (1-\tau^\alpha) \frac{\partial}{\partial t}(-\dot{e}_{\alpha,1}\ast_t \varrho \dot{\bu}) \\
&\qquad = \Div(2\mu \beps(\bu) + \lambda \tr(\beps(\bu)) \bI)\,
\\
&\qquad \quad +  (\tau^\alpha-1)\,  \dot{e}_{\alpha,1}\, \varrho\bh \,+\, e_{\alpha,1}\, \Div (\tau^\alpha\bs-
2\mu \beps(\bg) - \lambda \tr(\beps(\bg)) \bI)\\
&\qquad \quad +  \tau^\alpha \bF + (\tau^\alpha-1) \dot{e}_{\alpha,1} \ast_t \bF.
\end{aligned}
\end{align}
By introducing the function
\[ \bG:= (\tau^\alpha-1)\,  \dot{e}_{\alpha,1}\, \varrho \bh \,+\, e_{\alpha,1}\, \Div (\tau^\alpha\bs-
2\mu \beps(\bg) - \lambda \tr(\beps(\bg)) \bI) +  \tau^\alpha \bF + (\tau^\alpha-1) \dot{e}_{\alpha,1} \ast_t \bF\]
that collects the terms involving the initial data $\bg$, $\bh$, $\bs$ and the load vector $\bF$ on the right-hand side of \eqref{eq:7}, the equation \eqref{eq:7} takes the following more compact form:
\begin{align}\label{eq:8}
&\tau^\alpha \varrho \ddot \bu + (1-\tau^\alpha)\, \frac{\partial}{\partial t}(-\dot{e}_{\alpha,1}\ast_t \varrho \dot{\bu})(t) = \Div(2\mu \beps(\bu) + \lambda \tr(\beps(\bu)) \bI)\, + \, \bG.
\end{align}
We shall refer to  equation \eqref{eq:8} as the fractional Zener wave equation in three dimensional space.

 Next we shall derive a formal energy identity for the initial-boundary-value problem \eqref{eq:2}, \eqref{eq:3}, \eqref{eq:8}.

\section{Formal energy estimate for the model}
We begin the analysis of the problem by establishing a formal energy inequality, which we shall later rigorously prove by means of an abstract Galerkin approximation. We shall then use the energy inequality satisfied by the sequence of Galerkin approximations in conjunction with a compactness argument to show the existence of weak solutions to the initial-boundary-value problem \eqref{eq:2}, \eqref{eq:3}, \eqref{eq:8} under consideration, and we shall also prove the uniqueness of weak solutions. For the moment, though, we shall postulate the existence of
sufficiently smooth solutions in order to proceed with the formal derivation of an energy identity for the model.

To this end we shall take the scalar product of \eqref{eq:8} with $\dot{\bu}$, integrate the resulting equality over $\Omega$, and perform partial integration with respect to the spatial variable $\bx$, noting that $\bu$, and therefore also $\dot{\bu}$, satisfies a homogeneous Dirichlet boundary condition on $(0,T] \times \partial \Omega$.
In order to avoid notational clutter, whenever the function $\bF$ is extended by $\mathbf{0}$ from $(0,T] \times \Omega$ to
$(0,\infty) \times \Omega$  the extended function will be denoted by the same symbol as the original function.

As will be seen below, it is significant for the derivation of the energy identity, which guarantees continuous dependence of the solution on the data, that:
\begin{itemize}
\item $\tau \in (0,1]$, by hypothesis; and
\item $e_{\alpha,1}\geq 0$, $-\dot e_{\alpha,1} \geq 0$ and $\ddot e_{\alpha,1} \geq 0$ on $(0,T]$, with $\dot e_{\alpha,1} \in {\rm L}^1((0,T))$ and $\ddot{e}_{\alpha,1} \in {\rm L}^1_{\mathrm{loc}}((0,T))$ for all $T>0$.
\end{itemize}
We note in passing that by a similar reasoning the discussion below can be replicated in the case of the standard (integer-order) Zener model, corresponding to $\alpha=1$, but since the analysis of that model is much simpler we shall not include it here and will confine ourselves to the fractional-order Zener model, with $\alpha \in (0,1)$. An identical comment
applies to the case of a Hookean solid, corresponding to taking $\alpha =0$ in \eqref{eq:4}.

By formally testing the equation \eqref{eq:8} with $\dot{\bu}$ and noting that $\dot{\bu}$ satisfies a homogeneous Dirichlet
boundary condition on $(0,T] \times \partial \Omega$ we deduce, by partial integration with respect to the spatial variable
$\bx$, that, for any $t \in (0,T]$,
\begin{align*}
&\frac{\tau^\alpha}{2} \frac{\dd}{\dd t}\int_\Omega \varrho |\dot \bu(t,\bx)|^2 \dd \bx  + (1-\tau^\alpha)\, \int_\Omega \varrho  \frac{\partial}{\partial t}(-\dot{e}_{\alpha,1}\ast_t \dot{\bu})(t,\bx)\cdot \dot{\bu}(t,\bx) \dd \bx \\
&\qquad + \frac{1}{2} \frac{\dd}{\dd t} \int_\Omega 2\mu |\beps(\bu(t,\bx))|^2 + \lambda |\tr(\beps(\bu(t,\bx)))|^2  \dd \bx \, = \, \int_\Omega \bG(t,\bx) \cdot \dot{\bu}(t,\bx) \dd \bx.
\end{align*}
Hence, by integration over $t \in (0,T]$ and noting the initial conditions
\eqref{eq:2}, we deduce that
\begin{align}\label{eq:9}
\begin{aligned}
&\frac{\tau^\alpha }{2} \int_\Omega \varrho  |\dot \bu(t,\bx)|^2 \dd \bx  + (1-\tau^\alpha)\, \int_0^t \int_\Omega  \frac{\partial}{\partial s}(-\dot{e}_{\alpha,1}\ast_s \sqrt{\varrho} \dot{\bu})(s,\bx)\cdot \sqrt{\varrho} \dot{\bu}(s,\bx) \dd \bx \dd s \\
&\quad \quad + \frac{1}{2} \int_\Omega 2\mu |\beps(\bu(t,\bx))|^2 + \lambda |\tr(\beps(\bu(t,\bx)))|^2  \dd \bx \\
&= \, \int_0^t \int_\Omega \bG \cdot \dot{\bu} (s,\bx) \dd \bx \dd s + \frac{\tau^\alpha }{2} \int_\Omega \varrho |\bh(\bx)|^2 \dd \bx + \frac{1}{2} \int_\Omega 2\mu |\beps(\bg(\bx))|^2 + \lambda |\tr(\beps(\bg(\bx)))|^2  \dd \bx.
\end{aligned}
\end{align}
To proceed, we need to show that the second term on the left-hand side of \eqref{eq:9} is nonnegative, and that $\dot{\bu}$ can be eliminated from the right-hand side by absorbing it into the terms appearing on the left-hand side.
Once the nonnegativity of the second term on the left-hand side of \eqref{eq:9} has been verified, the identity \eqref{eq:9} can be viewed
as expressing balance of the total energy. In particular, when the load vector $\bF=\mathbf{0}$ and the initial data are such that $\bG=\mathbf{0}$, we have that
\begin{align}\label{eq:10}
\begin{aligned}
&\frac{\tau^\alpha }{2} \int_\Omega \varrho  |\dot \bu(t,\bx)|^2 \dd \bx  + (1-\tau^\alpha)\, \int_0^t \int_\Omega \frac{\partial}{\partial s}(-\dot{e}_{\alpha,1}\ast_s \sqrt{\varrho} \dot{\bu})(s,\bx)\cdot \sqrt{\varrho} \dot{\bu}(s,\bx) \dd \bx \dd s \\
&\qquad + \frac{1}{2}\int_\Omega 2\mu |\beps(\bu(t,\bx))|^2 + \lambda |\tr(\beps(\bu(t,\bx)))|^2  \dd \bx  \, \\
&= \frac{\tau^\alpha }{2} \int_\Omega \varrho \,|\bh(\bx)|^2 \dd \bx
+ \frac{1}{2} \int_\Omega 2\mu |\beps(\bg(\bx))|^2 + \lambda |\tr(\beps(\bg(\bx)))|^2  \dd \bx, \qquad t \in (0,T].
\end{aligned}
\end{align}
Even more specifically, if $\bF=\mathbf{0}$ and $\tau=1$, and $\bs$ is
related to $\beps(\bg)$  through Hooke's law (i.e., $\bs =  2\mu \beps(\bg) + \lambda \mathrm{tr}(\beps(\bg))$), whereby also $\bG= \mathbf{0}$, then the second term on the left-hand side of \eqref{eq:9} (which, thanks to Lemma \ref{le:1} below, can be viewed as an energy dissipation term,) is absent, as is the first term on the right-hand side of \eqref{eq:9}, and we have conservation of the total energy:
\[ \mathcal{E}(t):= \frac{1}{2} \int_\Omega \varrho |\dot \bu(t,\bx)|^2 \dd \bx
+ \frac{1}{2}\int_\Omega 2\mu |\beps(\bu(t,\bx))|^2 + \lambda |\tr(\beps(\bu(t,\bx)))|^2  \dd \bx   = \mathcal{E}(0) \quad \forall\,t \in [0,T].\]

Returning to the general case,
to show the nonnegativity of the second term on the left-hand side of \eqref{eq:9} we invoke the following result
(cf. Lemma 1.7.2 in \cite{Siskova}, whose proof is based on the identity stated in Lemma 2.3.1 in the work of Zacher \cite{Zacher}; see also identity (9) in \cite{Zacher0}).

\begin{lemma}\label{le:1}
Let $\mathcal{H}$ be a separable Hilbert space over the field of real numbers, with scalar product $(\cdot,\cdot)_{\mathcal{H}}$ and norm $\|\cdot\|_{\mathcal{H}}$, and let $T > 0$.
Then, for any $k \in \LL^1(0,T)$ such that
$k \geq 0$, $\dot{k} \in \LL^1_{\rm loc}(0,T)$, and $\dot{k} \leq 0$, and any
$v \in \LL^2((0, T);\mathcal{H})$, the following inequality holds:
\begin{align*}
\int_0^t \left(\frac{\dd}{\dd s}(k \ast_t v)(s), v(s)\right)_{\mathcal{H}} \dd s
& \geq  \frac{1}{2} (k \ast_t \|v(\cdot)\|^2_{\mathcal{H}})(t)
+ \frac{1}{2}\int_0^t k(s)\|v(s)\|^2_{\mathcal{H}} \dd s
\qquad \mbox{for all $t \in (0,T]$},
\end{align*}
each of the two terms on the right-hand side of the inequality being nonnegative.
\end{lemma}

Taking $k(t)=-\dot{e}_{\alpha,1}(t) (>0)$, $t \in (0,T]$, $\mathcal{H}=\LL^2_\varrho(\Omega)$, equipped
with the inner product and norm (and analogous notations for norms of weighted Lebesgue spaces, used in what follows, with weight functions $1/\varrho$, $\mu$, $1/\mu$, and $\lambda$ instead of $\varrho$) defined by
\[(\bv, \bw)_{\LL^2_\varrho(\Omega)}:=\int_\Omega \varrho(\bx)\, \bv(\bx) \cdot \bw(\bx)  \dd \bx, \qquad \|\bv\|_{L^2_\varrho(\Omega)}:=(\bv, \bv)^{\frac{1}{2}}_{\LL^2_\varrho(\Omega)},\]

\noindent
and $v = \dot{\bu}$ in Lemma \ref{le:1}, we deduce that the second term on the left-hand side of \eqref{eq:9}
is nonnegative.

It remains to show that the
function $\dot{\bu}$, appearing in the integrand of the first integral on the right-hand side, can
be absorbed into the left-hand side. To this end, we recall that
\[ \bG:= (\tau^\alpha-1)\,  \dot{e}_{\alpha,1}\, \varrho \bh \,+\, e_{\alpha,1}\, \Div (\tau^\alpha\bs-
2\mu \beps(\bg) - \lambda \tr(\beps(\bg)) \bI) +  \tau^\alpha \bF + (\tau^\alpha-1) \dot{e}_{\alpha,1} \ast_t \bF,\]
and we denote by $\bT_1$, $\bT_2$, $\bT_3$, and $\bT_4$, respectively, the four terms whose sum is $\bG$.

Clearly, because the function $t \in [0,\infty) \mapsto e_{\alpha,1}(t)$ is positive, strictly monotonic decreasing,
and $e_{\alpha,1}(0)=1$, we have by the Cauchy--Schwarz inequality that
\begin{align*}
\int_0^t \int_\Omega \bT_1(s,\bx) &\cdot \dot{\bu}(s,\bx) \dd \bx \dd s \leq  (1-\tau^\alpha)\,
\int_0^t (-\dot{e}_{\alpha,1}(s))
 \int_\Omega \varrho |\bh(\bx)| |\dot{\bu}(s,\bx)| \dd \bx \dd s\\
&\leq  (1-\tau^\alpha)\, \int_0^t (-\dot{e}_{\alpha,1}(s))\,
 \|\bh\|_{\LL^2_\varrho(\Omega)}\, \|\dot{\bu}(s,\cdot)\|_{\LL^2_\varrho(\Omega)} \dd s\\
&= (1-\tau^\alpha)\, \|\bh\|_{\LL^2_\varrho(\Omega)} \int_0^t (-\dot{e}_{\alpha,1}(s))\,
 \, \|\dot{\bu}(s,\cdot)\|_{\LL^2_\varrho(\Omega)} \dd s
\\
&
\leq (1-\tau^\alpha)\, \|\bh\|_{\LL^2_\varrho(\Omega)} \left(\int_0^t (-\dot{e}_{\alpha,1}(s)) \dd s\right)^{\frac{1}{2}}
\left(\int_0^t (-\dot{e}_{\alpha,1}(s)) \|\dot{\bu}(s,\cdot)\|^2_{\LL^2_\varrho(\Omega)} \dd s\right)^{\frac{1}{2}}\\
&
\leq (1-\tau^\alpha)\, \|\bh\|_{\LL^2_\varrho(\Omega)}
\left(e_{\alpha,1}(0) - e_{\alpha,1}(t)\right)^{\frac{1}{2}}
\left(\int_0^t (-\dot{e}_{\alpha,1}(s)) \|\dot{\bu}(s,\cdot)\|^2_{\LL^2_\varrho(\Omega)} \dd s\right)^{\frac{1}{2}}.
\end{align*}
By bounding the nonnegative factor $\left(e_{\alpha,1}(0) - e_{\alpha,1}(t)\right)^{\frac{1}{2}}$ above by $1$, for any $\delta_1>0$, to be fixed,
\begin{align}\label{eq:11}
\begin{aligned}
\int_0^t \int_\Omega \bT_1(s,\bx) &\cdot \dot{\bu}(s,\bx) \dd \bx \dd s \leq
(1-\tau^\alpha)\, \|\bh\|_{\LL^2_\varrho (\Omega)}
\left(\int_0^t (-\dot{e}_{\alpha,1}(s)) \|\dot{\bu}(s,\cdot)\|^2_{\LL^2_\varrho(\Omega)} \dd s\right)^{\frac{1}{2}}\\
&\leq
\frac{(1-\tau^\alpha)^2}{4\delta_1\tau^\alpha } \|\bh\|_{\LL^2_\varrho(\Omega)}^2 +
\tau^\alpha \delta_1 \int_0^t (-\dot{e}_{\alpha,1}(s)) \|\dot{\bu}(s,\cdot)\|^2_{\LL^2_\varrho(\Omega)} \dd s.
\end{aligned}
\end{align}

Next, by partial integration with respect to the temporal variable followed by partial integration
with respect to the spatial variable, we have, upon defining
\[ \bkap_0:= \tau^\alpha\bs-
2\mu \beps(\bg) - \lambda \tr(\beps(\bg)) \bI,\]
that
\begin{align*}
&\int_0^t \int_\Omega \bT_2(s,\bx) \cdot \dot{\bu}(s,\bx) \dd \bx \dd s = \int_0^t e_{\alpha,1}(s)\frac{\dd}{\dd s}\left[\int_\Omega \Div \bkap_0(\bx)\cdot \bu(s,\bx) \dd \bx\right] \dd s\\
&\qquad = \left[e_{\alpha,1}(s) \int_\Omega \Div \bkap_0(\bx) \cdot \bu(s,\bx) \dd \bx\right]_{s=0}^{s=t}
- \int_0^t \dot{e}_{\alpha,1}(s) \left[\int_\Omega \Div \bkap_0(\bx)\cdot \bu(s,\bx) \dd \bx\right] \dd s\\
&\qquad = \left[-e_{\alpha,1}(s) \int_\Omega \bkap_0(\bx) : \nabla \bu(s,\bx) \dd \bx\right]_{s=0}^{s=t}
+ \int_0^t \dot{e}_{\alpha,1}(s) \left[\int_\Omega \bkap_0(\bx) :  \nabla \bu(s,\bx) \dd \bx\right] \dd s.
\end{align*}
Now, letting $\mathbb{R}^{3 \times 3}_{\mathrm{sym}}$ denote the set of all symmetric $3 \times 3$
matrices with real entries, and noting that for any $A \in \mathbb{R}^{3 \times 3}_{\mathrm{sym}}$
and any $B \in \mathbb{R}^{3 \times 3}$ one has that $A:B = A: \frac{1}{2}(B + B^{\rm T})$, we deduce that
\begin{align*}
&\int_0^t \int_\Omega \bT_2(s,\bx) \cdot \dot{\bu}(s,\bx) \dd \bx \dd s
= \int_0^t \dot{e}_{\alpha,1}(s) \left[\int_\Omega \bkap_0(\bx) :  \beps(\bu(s,\bx)) \dd \bx\right] \dd s\\
& \qquad + \left[e_{\alpha,1}(0) \int_\Omega \bkap_0(\bx) : \beps(\bg(\bx)) \dd \bx\right]
- \left[e_{\alpha,1}(t) \int_\Omega \bkap_0(\bx) : \beps(\bu(t,\bx)) \dd \bx\right]
\\
&\leq \|\bkap_0\|_{\LL^2_{1/\mu}(\Omega)} \int_0^t (-\dot{e}_{\alpha,1}(s))  \|\beps(\bu(s,\cdot))\|_{\LL^2_\mu(\Omega)} \dd s
+  \|\bkap_0\|_{\LL^2_{1/\mu}(\Omega)}  \|\beps(\bg)\|_{\LL^2_\mu(\Omega)}\\
&\qquad + \|\bkap_0\|_{\LL^2_{1/\mu}(\Omega)} \|\beps(\bu(t))\|_{\LL^2_\mu(\Omega)},
\end{align*}
where in the transition to the right-hand side of the last inequality we have used that $e_{\alpha,1}(0)=1$ and that $t \in [0,\infty)
\mapsto e_{\alpha,1}(t)$ is positive and monotonic decreasing. Hence, by the Cauchy--Schwarz inequality,
and with a suitable real number $\delta_2>0$, to be fixed below,
\begin{align}\label{eq:12}
\begin{aligned}
&\int_0^t \int_\Omega \bT_2(s,\bx) \cdot \dot{\bu}(s,\bx) \dd \bx \dd s\\
&\quad\leq \|\bkap_0\|_{\LL^2_{1/\mu}(\Omega)} \left(\int_0^t (-\dot{e}_{\alpha,1}(s))\dd s\right)^{\frac{1}{2}}
\left(\int_0^t (-\dot{e}_{\alpha,1}(s))  \|\beps(\bu(s,\cdot))\|^2_{\LL^2_\mu(\Omega)}\dd s\right)^{\frac{1}{2}}
\\
& \qquad
+  \|\bkap_0\|_{\LL^2_{1/\mu}(\Omega)}  \|\beps(\bg)\|_{\LL^2_\mu(\Omega)}
+ \|\bkap_0\|_{\LL^2_{1/\mu}(\Omega)} \|\beps(\bu(t))\|_{\LL^2_\mu(\Omega)}\\
&\quad\leq \|\bkap_0\|_{\LL^2_{1/\mu}(\Omega)} \left(\int_0^t (-\dot{e}_{\alpha,1}(s))  \|\beps(\bu(s,\cdot))\|^2_{\LL^2_\mu(\Omega)}\dd s\right)^{\frac{1}{2}}
\\
& \qquad
+  \|\bkap_0\|_{\LL^2_{1/\mu}(\Omega)}  \|\beps(\bg)\|_{\LL^2_{\mu}(\Omega)}
+ \|\bkap_0\|_{\LL^2_{1/\mu}(\Omega)} \|\beps(\bu(t))\|_{\LL^2_\mu(\Omega)}\\
&\quad\leq \delta_2 \int_0^t (-\dot{e}_{\alpha,1}(s))  \|\beps(\bu(s,\cdot))\|^2_{\LL^2_\mu(\Omega)}\dd s
+ \frac{1}{4\delta_2}\|\bkap_0\|^2_{\LL^2_{1/\mu}(\Omega)}\\
&\qquad + \delta_2\|\beps(\bg)\|_{\LL^2_{\mu} (\Omega)}^2 + \frac{1}{4\delta_2}\|\bkap_0\|^2_{\LL^2_{1/\mu}(\Omega)} + \delta_2\|\beps(\bu(t))\|^2_{\LL^2_\mu(\Omega)} + \frac{1}{4\delta_2}\|\bkap_0\|^2_{\LL^2_{1/\mu}(\Omega)}.
\end{aligned}
\end{align}

Next, for a positive real number $\delta_{3}$, to be fixed below,
\begin{align}\label{eq:13}
\begin{aligned}
 \int_0^t \int_\Omega \bT_3(s,\bx) \cdot \dot{\bu}(s,\bx) \dd s \dd \bx &= \tau^\alpha \int_0^t \int_\Omega \bF(s,\bx) \cdot \dot{\bu}(s,\bx) \dd s \dd \bx
 \\& \leq  \tau^\alpha \int_0^t \|\bF(s)\|_{\LL^2_{1/\varrho}(\Omega)} \|\dot{\bu}(s)\|_{\LL^2_\varrho(\Omega)} \dd s
 \\& \leq \delta_3 \int_0^t \|\dot{\bu}(s)\|^2_{\LL^2_\varrho(\Omega)} \dd s
 + \frac{\tau^{2\alpha}}{4\delta_3}\int_0^t \|\bF(s)\|^2_{\LL^2_{1/\varrho}(\Omega)} \dd s.
\end{aligned}
\end{align}

Finally, by the Cauchy--Schwarz inequality with respect to $\bx$, Minkowski's integral inequality, the negativity of
$\dot{e}_{\alpha,1}$, the bound $\|-\dot{e}_{\alpha,1}\|_{\LL^1(0,t)} = 1 - e_{\alpha,1}(t) \leq 1$, Young's inequality for the (Laplace) convolution $\ast_t$ (whose proof we have
included at the end of this section for the sake of completeness; cf. Lemma \ref{le:2}), and with $\delta_4>0$
to be fixed below, we have that
\begin{align}\label{eq:14}
\begin{aligned}
\int_0^t \int_\Omega \bT_4(s,\bx) &\cdot \dot{\bu}(s,\bx) \dd s \dd \bx =(1-\tau^\alpha) \int_0^t \int_\Omega
(-\dot{e}_{\alpha,1} \ast_s \bF)(s,\bx) \cdot \dot{\bu}(s,\bx) \dd s \dd \bx\\
& \leq (1- \tau^\alpha) \int_0^t  \|-\dot{e}_{\alpha,1} \ast_s \bF(s)\|_{\LL^2_{1/\varrho}(\Omega)} \|\dot{\bu}(s)\|_{\LL^2_{\varrho}(\Omega)}
\dd s\\
&\leq (1- \tau^\alpha) \int_0^t  (-\dot{e}_{\alpha,1} \ast_s \|\bF\|_{\LL^2_{1/\varrho}(\Omega)})(s) \|\dot{\bu}(s)\|_{\LL^2_{\varrho}(\Omega)}
\dd s\\
& \leq \delta_4 \int_0^t \|\dot{\bu}(s)\|^2_{\LL^2_\varrho(\Omega)} \dd s + \frac{(1-\tau^\alpha)^2}{4\delta_4}\int_0^t
|(-\dot{e}_{\alpha,1} \ast_s \|\bF\|_{\LL^2_{1/\varrho}(\Omega)})(s)|^2 \dd s\\
& = \delta_4 \int_0^t \|\dot{\bu}(s)\|^2_{\LL^2_\varrho(\Omega)} \dd s + \frac{(1-\tau^\alpha)^2}{4\delta_4}
\left[\|(-\dot{e}_{\alpha,1}) \ast_s \|\bF\|_{\LL^2_{1/\varrho}(\Omega)}\|_{\LL^2(0,t)}\right]^2\\
& \leq  \delta_4 \int_0^t \|\dot{\bu}(s)\|^2_{\LL^2_\varrho(\Omega)} \dd s + \frac{(1-\tau^\alpha)^2}{4\delta_4}
\left[\|(-\dot{e}_{\alpha,1})\|_{\LL^1(0,t)} \|\|\bF\|_{\LL^2_{1/\varrho}(\Omega)}\|_{\LL^2(0,t)}\right]^2\\
& \leq  \delta_4 \int_0^t \|\dot{\bu}(s)\|^2_{\LL^2_\varrho(\Omega)} \dd s + \frac{(1-\tau^\alpha)^2}{4\delta_4}
\int_0^t \|\bF(s)\|^2_{\LL^2_{1/\varrho}(\Omega)}\dd s.
\end{aligned}
\end{align}
By substituting \eqref{eq:11}--\eqref{eq:14} into \eqref{eq:10} we deduce that
\begin{align}\label{eq:15}
\begin{aligned}
&\frac{\tau^\alpha}{2} \|\dot \bu(t)\|^2_{\LL^2_\varrho(\Omega)}  +
(1-\tau^\alpha) \int_0^t \int_\Omega
\frac{\partial}{\partial s}(-\dot{e}_{\alpha,1}\ast_s \sqrt{\varrho}\dot{\bu})(s,\bx)\cdot \sqrt{\varrho} \dot{\bu}(s,\bx) \dd s \dd \bx\\
&\qquad + \|\beps(\bu(t))\|^2_{\LL^2_\mu(\Omega)} + \frac{1}{2} \|\tr(\beps(\bu(t)))\|^2_{\LL^2_\lambda(\Omega)}\\
&\leq \frac{\tau^\alpha }{2} \|\bh\|^2_{\LL^2_\varrho(\Omega)}
+ \|\beps(\bg)\|^2_{\LL^2_\mu(\Omega)}  + \frac{1}{2} \|\tr(\beps(\bg))\|^2_{\LL^2_\lambda(\Omega)}\\
&\qquad + \tau^\alpha \delta_1 \int_0^t (-\dot{e}_{\alpha,1}(s)) \|\dot{\bu}(s)\|^2_{\LL^2_\varrho(\Omega)} \dd s
+ \frac{(1-\tau^\alpha)^2}{4\delta_1\tau^\alpha } \|\bh\|_{\LL^2_\varrho(\Omega)}^2\\
&\qquad + \delta_2 \int_0^t (-\dot{e}_{\alpha,1}(s))  \|\beps(\bu(s))\|^2_{\LL^2_\mu(\Omega)}\dd s
+ \frac{1}{4\delta_2}\|\bkap_0\|^2_{\LL^2_{1/\mu}(\Omega)}\\
&\qquad + \delta_2\|\beps(\bg)\|_{\LL^2_\mu(\Omega)}^2 + \frac{1}{4\delta_2}\|\bkap_0\|^2_{\LL^2_{1/\mu}(\Omega)}\\
& \qquad + \delta_2\|\beps(\bu(t))\|^2_{\LL^2_\mu(\Omega)} + \frac{1}{4\delta_2}\|\bkap_0\|^2_{\LL^2_{1/\mu}(\Omega)}\\
&\qquad + \delta_3 \int_0^t \|\dot{\bu}(s)\|^2_{\LL^2_\varrho(\Omega)} \dd s
 + \frac{\tau^{2\alpha}}{4\delta_3}\int_0^t \|\bF(s)\|^2_{\LL^2_{1/\varrho}(\Omega)} \dd s\\
&\qquad + \delta_4 \int_0^t \|\dot{\bu}(s)\|^2_{\LL^2_\varrho(\Omega)} \dd s + \frac{(1-\tau^\alpha)^2}{4\delta_4}
\int_0^t \|\bF(s)\|^2_{\LL^2_{1/\varrho}(\Omega)}\dd s.
\end{aligned}
\end{align}
We now fix
\[ \delta_1 = \delta_2 = \frac{1}{2},\quad
 \delta_3 = \delta_4 = \frac{\tau^\alpha}{4}.
\]
The inequality \eqref{eq:15} then takes the following form, for $t \in (0,T]$:
\begin{align}\label{eq:16}
\begin{aligned}
&\frac{\tau^\alpha}{2} \|\dot \bu(t)\|^2_{\LL^2_\varrho(\Omega)}  +
(1-\tau^\alpha) \int_0^t \int_\Omega
\frac{\partial}{\partial s}(-\dot{e}_{\alpha,1}\ast_s \sqrt{\varrho}\dot{\bu})(s,\bx)\cdot \sqrt{\varrho}\dot{\bu}(s,\bx) \dd s \dd \bx\\
&\qquad + \frac{1}{2} \|\beps(\bu(t))\|^2_{\LL^2_\mu(\Omega)} + \frac{1}{2} \|\tr(\beps(\bu(t)))\|^2_{\LL^2_\lambda(\Omega)}\\
&\leq \frac{\tau^{2\alpha} + (1-\tau^\alpha)^2}{2\tau^\alpha} \|\bh\|^2_{\LL^2_\varrho(\Omega)}
+ \frac{3}{2} \|\beps(\bg)\|^2_{\LL^2_\mu(\Omega)}  + \frac{1}{2} \|\tr(\beps(\bg))\|^2_{\LL^2_\lambda(\Omega)}
+ \frac{3}{2}\|\bkap_0\|^2_{\LL^2_{1/\mu}(\Omega)}\\
&\qquad + \frac{\tau^{2\alpha}+(1-\tau^\alpha)^2}{\tau^\alpha}\int_0^t \|\bF(s)\|^2_{\LL^2_{1/\varrho}(\Omega)} \dd s\\
&\qquad + \frac{\tau^\alpha}{2} \int_0^t (1-\dot{e}_{\alpha,1}(s)) \|\dot{\bu}(s)\|^2_{\LL^2_\varrho(\Omega)} \dd s
+ \frac{1}{2} \int_0^t (-\dot{e}_{\alpha,1}(s))  \|\beps(\bu(s))\|^2_{\LL^2_\mu(\Omega)}\dd s.
\end{aligned}
\end{align}
Now, consider the following two nonnegative functions defined on $[0,T]$:
\begin{align*}
y(t)&:= \frac{\tau^\alpha}{2} \|\dot \bu(t)\|^2_{\LL^2_\varrho(\Omega)} + \frac{1}{2} \|\beps(\bu(t))\|^2_{\LL^2_\mu(\Omega)} + \frac{1}{2} \|\tr(\beps(\bu(t)))\|^2_{\LL^2_\lambda(\Omega)},\\
z(t)&:= (1-\tau^\alpha) \int_0^t \int_\Omega
\frac{\partial}{\partial s}(-\dot{e}_{\alpha,1}\ast_s \sqrt{\varrho}\dot{\bu})(s,\bx)\cdot \sqrt{\varrho}\dot{\bu}(s,\bx) \dd s \dd \bx,
\end{align*}
and let
\begin{align}\label{eq:17}
\begin{aligned}
A(t) &:= \frac{\tau^{2\alpha} + (1-\tau^\alpha)^2}{2\tau^\alpha} \|\bh\|^2_{\LL^2_\varrho(\Omega)}
+ \frac{3}{2} \|\beps(\bg)\|^2_{\LL^2_\mu(\Omega)}  + \frac{1}{2} \|\tr(\beps(\bg))\|^2_{\LL^2_\lambda(\Omega)}
\\
&\qquad + \frac{3}{2}\|\bkap_0\|^2_{\LL^2_{1/\mu}(\Omega)} + \frac{\tau^{2\alpha}+(1-\tau^\alpha)^2}{\tau^\alpha}\int_0^t \|\bF(s)\|^2_{\LL^2_{1/\varrho}(\Omega)} \dd s.
\end{aligned}
\end{align}
Clearly, $0 \leq A(t) \leq A(T)=:A$. The inequality \eqref{eq:16} then implies that
\[ y(t) + z(t) \leq A(t) + \int_0^t (1-\dot{e}_{\alpha,1}(s)) y(s) \dd s.\]
Since $t \in [0,T] \mapsto A(t)$ is a nonnegative and nondecreasing function, by Gronwall's lemma we have that
\[ y(t) + z(t) \leq A(t)\, \mathrm{exp}\left(\int_0^t (1-\dot{e}_{\alpha,1}(s))\dd s\right) = A(t)\, \mathrm{exp}(t + 1 - e_{\alpha,1}(t)),
\quad t \in (0,T].\]
In other words,  with
\[A(t)=A(\tau^\alpha,\|\bh\|_{\LL^2_\varrho(\Omega)},\|\beps(\bg)\|_{\LL^2_\mu(\Omega)},
\|\tr(\beps(\bg))\|_{\LL^2_\lambda(\Omega)}, \|\bkap_0\|_{\LL^2_{1/\mu}(\Omega)}, \|\bF\|_{\LL^2(0,t;\LL^2_{1/\varrho}(\Omega))})\geq 0\]
defined by the expression \eqref{eq:17} for $t \in [0,T]$, the following energy inequality holds for all $t \in [0,T]$:
\begin{align}\label{eq:18}
\begin{aligned}
&\frac{\tau^\alpha}{2} \|\dot \bu(t)\|^2_{\LL^2_\varrho(\Omega)} + \frac{1}{2} \|\beps(\bu(t))\|^2_{\LL^2_\mu(\Omega)} + \frac{1}{2} \|\tr(\beps(\bu(t)))\|^2_{\LL^2_\lambda(\Omega)},\\
&\quad + (1-\tau^\alpha) \int_0^t \int_\Omega
\frac{\partial}{\partial s}(-\dot{e}_{\alpha,1}\ast_s \sqrt{\varrho} \dot{\bu})(s,\bx)\cdot \sqrt{\varrho} \dot{\bu}(s,\bx) \dd s \dd \bx
\leq A(t)\, \mathrm{exp}(t+1 - e_{\alpha,1}(t)).
\end{aligned}
\end{align}
Thus, assuming the existence of a (sufficiently smooth) solution
$\bu$ to \eqref{eq:2}, \eqref{eq:3}, \eqref{eq:8}, with
\begin{align}\label{eq:19}
\bg \in [\HH^1_0(\Omega)]^3, \quad \bh \in [\LL^2(\Omega)]^3, \quad \bs = \bs^{\rm T}
\in [\LL^2(\Omega)]^{3 \times 3}, \quad \bF \in \LL^2(0,T;[\LL^2(\Omega)]^3),
\end{align}
recalling that, by hypothesis \textcolor{black}{\eqref{coeff-ass}}, $\varrho, \mu, \lambda \in \LL^\infty(\Omega)$, $\varrho$ and $\mu$ are bounded below by positive
constants $\varrho_0$ and $\mu_0$, respectively, and $\lambda \geq 0$ a.e. on $\Omega$,
the energy inequality \eqref{eq:18} holds, with $A(t)<\infty$ for all $t \in [0,T]$.

We emphasize here the significance of our assumption that $\tau \in (0,1]$: the positivity of $\tau$ is necessary in order to ensure that the factor $A(t)$ (cf. \eqref{eq:17}) appearing on the right-hand side of the energy inequality \eqref{eq:18} is finite, while $\tau \leq 1 ~(=\rho)$ ensures that the prefactor of the last term on the left-hand side of \eqref{eq:18}, which can be viewed as a nonnegative energy dissipation term thanks to Lemma \ref{le:1}, is nonnegative, whereby the entire left-hand side of \eqref{eq:18} is nonnegative.

\begin{remark}\label{re:1}
We remark that if $\bs$ is chosen so that $\tau^\alpha \bs = 2 \mu \beps(\bg)+ \lambda \tr(\beps(\bg))$,
then $\bkap_0 = \mathbf{0}$, and therefore also $\bT_2=\mathbf{0}$. The energy inequality
\eqref{eq:18} is then simpler and sharper, which can be seen by erasing all terms containing $\delta_2$
from the right-hand side of \eqref{eq:15}, and making the same choices of $\delta_1$, $\delta_3$ and $\delta_4$ as above.
{\color{black} In the special case of $\lambda=0$ this particular choice of the initial stress $\bs$, namely $\bs = 2\mu (1/\tau)^\alpha \beps(\bg)$, in our initial condition \eqref{eq:2}$_3$ results in the same initial condition as the one stated in equation (13) in the work of Freed and Diethelm \cite{FS} (recall that we scaled $\rho$ to $1$, so $(\rho/\tau)^\alpha = (1/\tau)^\alpha$).}
We shall proceed without making this restrictive assumption on $\bs$, and continue to study the general case when
$\tau^\alpha \bs$ is not required to be equal to
$2 \mu \beps(\bg)+ \lambda \tr(\beps(\bg))$.
\end{remark}

In the next section we shall use a compactness argument, based on a sequence of spatial Galerkin approximations
to the problem, to show the existence of a (unique) weak solution.

We close this section with the proof of Young's inequality for Laplace-type
convolution, which we used in the derivation of the energy inequality. The proof of this result in the case of Fourier-type convolution is standard;
in the case of Laplace-type convolution the argument proceeds along similar
lines, with minor modifications; {\color{black} we have included its statement and proof for the convenience of the reader.}

\begin{lemma}\label{le:2}
Let {\color{black} $p, q, r \in [1,\infty]$} be such that $\frac{1}{p} + \frac{1}{q} - 1 = \frac{1}{r}$, and
let $f \in \LL^p(0,t)$ and $g \in \LL^q(0,t)$ for some $t>0$; then $s \in [0,t] \mapsto
(f\ast_s g)(s):=\int_0^s f(s-u) g(u) \dd u \in \LL^r(0,t)$, and
\[ \|f \ast_s g \|_{\LL^r(0,t)} \leq \|f\|_{\LL^p(0,t)} \|g\|_{\LL^q(0,t)}.\]
\end{lemma}

\noindent\textit{Proof.}
{\color{black} If $p=\infty$, then necessarily $q=1$ and $r=\infty$, and if $q=\infty$, then necessarily $p=1$ and $r=\infty$.
Since for $r=\infty$ the result is a direct consequence of H\"older's
inequality,} we shall concentrate here on the nontrivial case
when $p,q,r \in [1,\infty)$. We begin by noting that because
\[ \frac{1}{r} + \frac{r-p}{pr} + \frac{r-q}{qr} = 1,\]
we have by H\"older's inequality that,  for any $s \in (0,t]$,
\begin{align*}
|(f \ast_s g)(s)| &= \left|\int_0^s f(s-u) g(u) \dd u\right| \leq \int_0^s |f(s-u)|\, |g(u)| \dd u \\
& = \int_0^s |f(s-u)|^{\frac{p}{r}} |g(u)|^{\frac{q}{r}}
|f(s-u)|^{1-\frac{p}{r}} |g(u)|^{1-\frac{q}{r}} \dd u\\
& \leq \||f(s-\cdot)|^{\frac{p}{r}} |g(\cdot)|^{\frac{q}{r}}\|_{\LL^r(0,s)}
\||f(s-\cdot)|^{1-\frac{p}{r}}\|_{\LL^{\frac{pr}{r-p}}(0,s)} \||g(\cdot)|^{1-\frac{q}{r}}\|_{\LL^{\frac{qr}{r-q}}(0,s)}
\\
& = \left(\int_0^s |f(s-u)|^{p} |g(u)|^{q} \dd u\right)^{\frac{1}{r}}
\|f\|^{\frac{r-p}{r}}_{\LL^p(0,s)} \|g\|^{\frac{r-q}{r}}_{\LL^q(0,s)}\\
& \leq \left(\int_0^s |f(s-u)|^{p} |g(u)|^{q} \dd u\right)^{\frac{1}{r}}
\|f\|^{\frac{r-p}{r}}_{\LL^p(0,t)} \|g\|^{\frac{r-q}{r}}_{\LL^q(0,t)}.
\end{align*}
Hence, by integration over $s \in (0,t)$, applying Fubini's theorem, and performing
the change of variable $\sigma := s-u$, we deduce that
\begin{align*}
\int_0^t |(f \ast_s g)(s)|^r \dd s
& \leq  \left(\int_0^t \int_0^s |f(s-u)|^{p} |g(u)|^{q} \dd u \dd s\right)
\|f\|^{r-p}_{\LL^p(0,t)} \|g\|^{r-q}_{\LL^q(0,t)}\\
& =   \left(\int_0^t |g(u)|^{q} \left(\int_u^t |f(s-u)|^{p}  \dd s\right) \dd u\right)
\|f\|^{r-p}_{\LL^p(0,t)} \|g\|^{r-q}_{\LL^q(0,t)}\\
& =   \left(\int_0^t |g(u)|^{q} \left(\int_0^{t-u} |f(\sigma)|^{p}  \dd \sigma\right) \dd u\right)
\|f\|^{r-p}_{\LL^p(0,t)} \|g\|^{r-q}_{\LL^q(0,t)}\\
& \leq \left(\int_0^t |g(u)|^{q} \left(\int_0^{t} |f(\sigma)|^{p}  \dd \sigma\right) \dd u\right)
\|f\|^{r-p}_{\LL^p(0,t)} \|g\|^{r-q}_{\LL^q(0,t)}\\
& = \left(\|g\|^{q}_{\LL^q(0,t)} \|f\|^{p}_{\LL^p(0,t)}\right)
\|f\|^{r-p}_{\LL^p(0,t)} \|g\|^{r-q}_{\LL^q(0,t)} = \|f\|^{r}_{\LL^p(0,t)} \|g\|^{r}_{\LL^q(0,t)}.
\end{align*}
By raising this to the power $\frac{1}{r}$, we
arrive at the desired inequality. $\quad\Box$

\section{Existence of weak solutions}
 Hereafter $\WW^{s,p}(D)$ will denote the Sobolev space of real-valued functions defined on a bounded open set $D \subset \mathbb{R}^d$, $d \geq 1$,  with differentiability index $s>0$ and integrability index $p \in [1,\infty]$ (cf. \cite{AF}). When $p=2$, we shall write $\HH^s(\Omega)$ instead of $\WW^{s,2}(D)$ and $\HH^s_0(D)$ will denote the closure of $\CC^\infty_0(D)$ in $\HH^s(D)$. When $D$ is a bounded open Lipschitz domain and $s \in (\frac{1}{2},\frac{3}{2})$, elements of $\HH^s_0(D)$ have
zero trace on $\partial D$; for such $s$, $\HH^{-s}(D)$ will denote the dual space of $\HH^s_0(D)$.

For a Banach
space $\mathcal{B}$, we shall denote by $\LL^p(0,T;\mathcal{B})$ and $\WW^{s,p}(0,T;\mathcal{B})$, respectively, the associated Lebesgue and Sobolev space of $\mathcal{B}$-valued mappings defined on the open interval $(0,T)$, and $\CC([0,T];\mathcal{B})$ will signify the set of all uniformly continuous $\mathcal{B}$-valued functions defined on $[0,T]$. Furthermore, $\CC^{0,1}([0,T];\mathcal{B})$ will denote the space of Lipschitz-continuous $\mathcal{B}$-valued functions defined on $[0,T]$.
Suppose that $\mathcal{H}$ is a Hilbert space over the field of real numbers with
inner product $(\cdot,\cdot)_{\mathcal{H}}$. We shall denote by $\CC_w([0,T];\mathcal{H})$ the linear space of all weakly continuous functions from $[0,T]$ into $\mathcal{H}$, i.e., the set of all
functions $v \in \LL^\infty(0,T;\mathcal{H})$ such that $t \in [0,T] \mapsto (v(t),w) \in
\mathbb{R}$ is a continuous function on $[0,T]$ for each $w \in \mathcal{H}$.

Our objective in this section is to show the existence and the uniqueness of a \textit{weak solution} to the problem  \eqref{eq:2}, \eqref{eq:3}, \eqref{eq:8}, defined as follows.

\begin{definition}[Weak solution]\label{weak-sol}
Suppose that the initial data $\bg$, $\bh$, $\bs$ and the source term $\bF$ satisfy \eqref{eq:19}, and assume that
$\tau \in (0,1]$,  $\alpha \in (0,1)$, and $\varrho$, $\mu$, and $\lambda$ are as in \eqref{coeff-ass}. A function
\begin{alignat}{2}\label{eq:20}
\begin{aligned}
\bu &\in \CC_w([0,T];[\HH^1_0(\Omega)]^3),\qquad \mbox{with}\\
\dot{\bu} \in \CC_w([0,T]; &~ [\LL^2(\Omega)]^3), \quad \mbox{and} \quad
(-\dot e_{\alpha,1})^{\frac{1}{2}}~\! \dot{\bu} \in \LL^2(0,T;[\LL^2(\Omega)]^3),
\end{aligned}
\end{alignat}
satisfying the equality
\begin{align}\label{eq:21}
\begin{aligned}
 \tau^\alpha \int_0^T &(\varrho \bu (s,\cdot), \ddot{\bv}(s,\cdot)) \dd s - (1-\tau^\alpha) \int_0^T ((-\dot{e}_{\alpha,1}\ast_s \varrho \dot{\bu})(s,\cdot), \dot{\bv}(s,\cdot)) \dd s \\
&\quad + \int_0^T \big ( 2\mu \beps(\bu(s,\cdot)) + \lambda \tr(\beps(\bu(s,\cdot))) \bI
\,,  \beps(\bv(s,\cdot)) \big) \dd s \\
& = - \tau^\alpha  (\varrho \bg, \dot{\bv}(0,\cdot)) + \tau^\alpha  (\varrho \bh, \bv(0,\cdot)) + \int_0^T \langle \bG(s,\cdot), \bv(s,\cdot) \rangle \dd s
\end{aligned}
\end{align}
for all $\bv \in \WW^{2,1}(0,T;[\LL^2(\Omega)]^3) \cap \LL^1(0,T;[\HH^1_0(\Omega)]^3)$
with $\bv(T,\cdot)=0$ and $\dot{\bv}(T,\cdot)=0$, and
\begin{align}\label{eq:22}
\bG:= (\tau^\alpha-1)\,  \dot{e}_{\alpha,1}\,\varrho \bh \,+\, e_{\alpha,1}\, \Div (\tau^\alpha\bs-
2\mu \beps(\bg) - \lambda \tr(\beps(\bg)) \bI) +  \tau^\alpha \bF + (\tau^\alpha-1) \dot{e}_{\alpha,1} \ast_t \bF,
\end{align}
is called a weak solution to the problem \eqref{eq:2}, \eqref{eq:3}, \eqref{eq:8}.
\end{definition}

In \eqref{eq:21} and throughout the rest of the paper $\langle \cdot , \cdot \rangle$ denotes the duality pairing between $[\HH^{-1}(\Omega)]^3$ and $[\HH^1_0(\Omega)]^3$,
and $(\cdot, \cdot)$ is the inner product of $[\LL^2(\Omega)]^3$. We note that, for $\alpha \in (0,1)$,
\begin{align}\label{eq:23}
 - \dot{e}_{\alpha,1}(t) \thicksim \frac{\alpha\, t^{\alpha-1}}{\Gamma(\alpha +1)}\qquad \mbox{as
 $t \rightarrow 0_+$},
\end{align}
and hence, by noting from \eqref{eq:22} the additive structure of $\bG$, we have that
\[\bG \in \LL^p(0,T;[\LL^2(\Omega)]^3) + \WW^{1,p}(0,T;[\HH^{-1}(\Omega)]^3) +
\LL^2(0,T;[\LL^2(\Omega)]^3)\qquad \forall\,
p \in \big[1, \textstyle{\frac{1}{1-\alpha}}\big).\]

{\color{black}
The function $\bsig$ has been eliminated in the transition from \eqref{eq:2}, \eqref{eq:3}, \eqref{eq:8}
to the weak formulation \eqref{eq:21}, \eqref{eq:22}, and
the initial condition $\bsig(0,\cdot) = \bs(\cdot)$ has been encoded into \eqref{eq:21}, \eqref{eq:22}.
Motivated by \eqref{eq:sig-laplace}, for a weak solution $\bu$, whose existence and uniqueness we will
show in Theorem \ref{th:1} below, we therefore \emph{define} the associated stress
tensor $\bsig$ by
\begin{align}\label{eq:stress-defin}
\begin{aligned}
\bsig(t,\cdot) &:= \mathcal{L}^{-1}\left(\frac{1+ p^\alpha}{1  + \tau^\alpha p^\alpha}\right) \ast_t
(2\mu \beps(\bu(t,\cdot)) + \lambda \tr(\beps(\bu(t,\cdot))) \bI) \\
&\quad + \mathcal{L}^{-1}\left(\frac{p^{\alpha-1}}{1  + \tau^\alpha p^\alpha}\right) \,(\tau^\alpha\bs(\cdot)-
2\mu \beps(\bg(\cdot)) - \lambda \tr(\beps(\bg(\cdot))) \bI).
\end{aligned}
\end{align}
}

Consider the bilinear form $a(\cdot,\cdot)$ on $[\HH^1_0(\Omega)]^3 \times [\HH^1_0(\Omega)]^3$, defined by
\[ a(\bw, \bv):= (2\mu \beps(\bw) + \lambda \tr(\beps(\bw)) \bI, \beps(\bv)) \qquad \forall\, \bw, \bv \in [\HH^1_0(\Omega)]^3,\]
and observe that
\[ a(\bw, \bv)= (2\mu \beps(\bw), \beps(\bv)) + (\lambda \tr(\beps(\bw)) , \tr(\beps(\bv))) \qquad \forall\, \bw, \bv \in [\HH^1_0(\Omega)]^3.\]
Clearly, $a(\bw,\bv) =  a(\bv,\bw)$, and there exist positive real numbers $c_1$ and $c_0$ such that $a(\bw,\bv)\leq c_1
\|\bw\|_{\HH^1(\Omega)} \|\bv\|_{\HH^1(\Omega)}$ for all $\bw, \bv \in [\HH^1_0(\Omega)]^3$ (by the Cauchy--Schwarz inequality), and $a(\bv,\bv) \geq c_0 \|\bv\|^2_{\HH^1(\Omega)}$ for all $\bv \in [\HH^1_0(\Omega)]^3$  (by Korn's inequality).
Hence, $a(\cdot,\cdot)$ is a symmetric, bounded, and coercive bilinear form on $[\HH^1_0(\Omega)]^3 \times [\HH^1_0(\Omega)]^3$. Furthermore, by Rellich's theorem, the infinite-dimensional separable Hilbert space
$[\HH^1_0(\Omega)]^3$ is compactly and densely embedded into the infinite-dimensional separable Hilbert space
$[\LL^2(\Omega)]^3$.

To proceed, we require the following version of the Hilbert--Schmidt theorem \cite{FS}.
\begin{lemma}\label{le:3} Let $\mathcal{H}$ and $\mathcal{V}$ be separable
Hilbert spaces, with $\mathcal{V}$ compactly embedded into $\mathcal{H}$ and $\overline{\mathcal{V}} = \mathcal{H}$
in the norm of $\mathcal{H}$. Let $a\colon \mathcal{V} \times \mathcal{V} \to \mathbb{R}$ be a nonzero, symmetric,
bounded and coercive bilinear form. Then, there exist sequences of real numbers
$(\lambda_n)_{n \in \mathbb{N}}$ and unit $\mathcal{H}$-norm members $(e_n)_{n \in \mathbb{N}}$ of $\mathcal{V}$,
which solve the following problem: \emph{Find $\lambda \in \mathbb{R}$ and
$e \in \mathcal{H} \setminus \{ 0 \}$ such that}
\begin{equation}\label{eq:24}
a(e,v) = \lambda ( e, v )_\mathcal{H} \quad \forall\,v \in \mathcal{V}.
\end{equation}
The $\lambda_n$, which can be assumed to be in increasing order with respect to $n$,
are positive, bounded from below away from $0$, and
$\lim_{n\to\infty}\lambda_n = \infty$.

Additionally, the $e_n$ form an $\mathcal{H}$-orthonormal system whose $\mathcal{H}$-closed span is
$\mathcal{H}$ and the rescaling $e_n/\sqrt{\lambda_n}$ gives rise to an $a$-orthonormal
system whose $a$-closed span is $\mathcal{V}$.
\end{lemma}

We are now ready to formulate the main result of this section.
\begin{theorem}\label{th:1}
Suppose that the initial data $\bg$, $\bh$, $\bs$ and the source term $\bF$ satisfy \eqref{eq:19}, and assume that
$\tau \in (0,1]$,  $\alpha \in (0,1)$, and $\varrho$, $\mu$, and $\lambda$ are as in \eqref{coeff-ass}.
Then, {\color{black} the weak formulation \eqref{eq:21}, \eqref{eq:22} of the
problem \eqref{eq:2}, \eqref{eq:3}, \eqref{eq:8} has a (weak) solution} in the sense of Definition \ref{weak-sol}
such that
\[\bu \in \CC([0,T];[\HH^s_0(\Omega)]^3)\qquad \mbox{for all $s \in (\frac{1}{2},1)$},\]
and
\[\tau^\alpha\varrho \dot \bu + (1-\tau^\alpha)(-\dot{e}_{\alpha,1}\ast_t \varrho \dot{\bu}) \in \WW^{1,p}(0,T;[\HH^{-1}(\Omega)]^3),\qquad \alpha \in (0,1),\]
for all $p \in [1,2]$ satisfying $p<\frac{1}{1-\alpha}$. Furthermore, $\bu$ satisfies the energy inequality
\begin{align}\label{EE}
\begin{aligned}
&\frac{\tau^\alpha}{2} \|\dot \bu(t')\|^2_{\LL^2_\varrho(\Omega)} + \frac{1}{2} \|\beps(\bu(t'))\|^2_{\LL^2_\mu(\Omega)} + \frac{1}{2} \|\tr(\beps(\bu(t')))\|^2_{\LL^2_\lambda(\Omega)}\\
&\quad + \frac{1-\tau^\alpha}{2} \int_0^{t'} -\dot{e}_{\alpha,1}(s)\|\dot{\bu}(s)\|^2_{\LL^2_\varrho(\Omega)} \dd s
\leq 3A(t) \, \mathrm{exp}(t+1 - e_{\alpha,1}(t)),
\end{aligned}
\end{align}
for all $t \in (0,T]$ and a.e. $t' \in (0,t]$, where $A(t)$ is defined by \eqref{eq:17} for $t \in [0,T]$.

{\color{black} The initial condition $\bu(0,\cdot) = \bg(\cdot)$ is satisfied in the sense of continuous functions
from $[0,T]$ into $[\LL^2(\Omega)]^3$ and the initial condition $\dot\bu (0,\cdot) = \bh(\cdot)$ is satisfied
as an equality in $\CC_w([0,T],[\LL^2(\Omega)]^3)$. Furthermore, the weak solution $\bu$ is unique and depends continuously on the data $\bg$, $\bh$, $\bs$, and $\bF$.

The stress tensor $\bsig$, defined by \eqref{eq:stress-defin} in terms of the unique weak solution $\bu$ of \eqref{eq:21}, \eqref{eq:22}, satisfies the initial condition $\bsig(0,\cdot) = \bs(\cdot)$ as an equality in $\CC_w([0,T],[\LL^2(\Omega)]^{3 \times 3})$.}
\end{theorem}

\textit{Proof.} STEP 1: \textit{Existence of solutions.} We begin by showing the existence of a weak solution.
We shall use Lemma \ref{le:3} with $\mathcal{H}= [\LL^2_\varrho(\Omega)]^3 \simeq [\LL^2(\Omega)]^3$ equipped
with the inner product defined by $(\bw,\bv)_{\mathcal{H}} := (\varrho \bw , \bv)$,  $\mathcal{V}= [\HH^1_0(\Omega)]^3$, to generate an $\mathcal{H}$-orthonormal Galerkin basis $(\mathbf{\bphi}_n)_{n \in \mathbb{N}}
\subset [\HH^1_0(\Omega)]^3$, whose $[\LL^2(\Omega)]^3$-closed span is
$[\LL^2(\Omega)]^3$ and the rescaling $\bphi_n/\sqrt{\lambda_n}$ gives rise to an $a$-orthonormal
system whose $a$-closed span is $[\HH^1_0(\Omega)]^3$; $(\lambda_n)_{n \in \mathbb{N}}$ is
a countably infinite sequence of positive eigenvalues, bounded away from $0$, and
$\lim_{n\to\infty}\lambda_n = \infty$, defined by $a(\bphi_n, \bv) = \lambda_n (\varrho \bphi_n , \bv)$
for all $\bv \in [\HH^1_0(\Omega)]^3$.

Let $\mathcal{V}_n := \mbox{span}\{\bphi_1,\dots, \bphi_n\}$, and let $P_n \bv \in \mathcal{V}_n$ denote
the orthogonal projection of $\bv \in [\LL^2_\varrho(\Omega)]^3$, in the inner product of $[\LL^2_\varrho(\Omega)]^3$, onto $\mathcal{V}_n$.  We seek a Galerkin approximation $\bu_n: [0,T] \mapsto
\bu_n(t) \in \mathcal{V}_n$ of the form
\begin{align}\label{eq:25}
\bu_n(t,\bx) := \sum_{k=1}^n \beta_k(t) \bphi_k(\bx)
\end{align}
satisfying
\begin{align}\label{eq:26}
\hspace{-3mm}\tau^\alpha (\varrho \ddot \bu_n, \bv ) + (1-\tau^\alpha) \left ( \frac{\partial}{\partial t}(-\dot{e}_{\alpha,1}\ast_t \varrho \dot{\bu}_n), \bv \right ) + \big ( 2\mu \beps(\bu_n) + \lambda \tr(\beps(\bu_n)) \bI
\,, \beps(\bv) \big) = \langle \bG, \bv \rangle
\end{align}
for all $\bv \in \mathcal{V}_n$, together with the initial conditions
\[ \bu_n(0,\cdot) = P_n\bg \qquad \mbox{and} \qquad \dot{\bu}_n(0,\cdot) = P_n\bh.\]
Equivalently,
\[\beta_k(0)=(\varrho \bg,\bphi_k)\qquad \mbox{and}\qquad \dot{\beta}_k(0)=(\varrho \bh,\bphi_k),\qquad
\mbox{for $k=1,\dots,n$}.\]

Hence,
\[ \|{\bu}_n(0,\cdot)\|_{\LL^2_\varrho(\Omega)} \leq \|\bg\|_{\LL^2_\varrho(\Omega)}\qquad \mbox{and}\qquad \|\dot{\bu}_n(0,\cdot)\|_{\LL^2_\varrho(\Omega)} \leq \|\bh\|_{\LL^2_\varrho(\Omega)},\]
and
\begin{align*}
a(\bu_n(0,\cdot),\bu_n(0,\cdot)) &= \sum_{k,\ell=1}^n \beta_k(0) \beta_\ell(0) a(\bphi_k,\bphi_\ell)
= \sum_{k,\ell=1}^n \beta_k(0) \beta_\ell(0) \lambda_k (\varrho  \bphi_k,\bphi_\ell)\\
&= \sum_{k=1}^n [\beta_k(0)]^2 \lambda_k \|\bphi_k\|^2_{\LL^2_\varrho(\Omega)}= \sum_{k=1}^n [(\varrho\bg,\bphi_k)]^2 \lambda_k \|\bphi_k\|^2_{\LL^2_\varrho(\Omega)}\\
&\leq \sum_{k=1}^\infty [(\varrho \bg,\bphi_k)]^2 \lambda_k \|\bphi_k\|^2_{\LL^2_\varrho(\Omega)}
= \sum_{k,\ell=1}^\infty (\varrho \bg,\bphi_k) (\varrho \bg, \bphi_\ell) \lambda_k
(\varrho \bphi_k,\bphi_\ell)\\
& = a\left(\sum_{k=1}^\infty (\varrho \bg,\bphi_k) \bphi_k , \sum_{\ell=1}^\infty (\varrho \bg,\bphi_\ell)
\bphi_\ell\right) = a(\bg,\bg).
\end{align*}
Thus, by the coercivity and the boundedness of the bilinear form $a(\cdot,\cdot)$ on
$[\HH^1_0(\Omega)]^3 \times [\HH^1_0(\Omega)]^3$, also
\[ c_0 \|\bu_n(0,\cdot)\|^2_{\HH^1(\Omega)} = c_0 \|P_n \bg\|^2_{\HH^1(\Omega)}\leq c_1 \|\bg\|^2_{\HH^1(\Omega)}.\]
Therefore, the orthogonal projector $P_n$ has operator norm $\|P_n\|_{\mathcal{L}([\LL^2_\varrho(\Omega)]^3, [\LL^2_\varrho(\Omega)]^3)}$
bounded by $1$, uniformly in $n$, and it is, simultaneously, a bounded linear operator from $[\HH^1_0(\Omega)]^3$ into $\mathcal{V}_n \subset [\HH^1_0(\Omega)]^3$, with operator norm $\|P_n\|_{\mathcal{L}([\HH^1(\Omega)]^3, [\HH^1(\Omega)]^3)}$ bounded by $(c_1/c_0)^{1/2}$, uniformly in $n$.

We begin by showing the existence of a unique Galerkin approximation $t \in [0,T] \mapsto \bu_n(t) \in \mathcal{V}_n$. By substituting \eqref{eq:25} into \eqref{eq:26} and taking $\bv = \bphi_m \in \mathcal{V}_n$ for $m=1,\dots, n$ and noting
the orthonormality $(\varrho \bphi_k, \bphi_m) = \delta_{k,m}$ for $k, m = 1, \dots, n$, we have that
\begin{align}\label{eq:27}
\tau^\alpha \ddot \beta_m  + (1-\tau^\alpha) \frac{\dd}{\dd t}(-\dot{e}_{\alpha,1}\ast_t \dot{\beta}_m)
+ \lambda_m \beta_m = \langle \bG, \bphi_m \rangle,\qquad \mbox{$m=1,\dots,n$,}
\end{align}
with $\langle \bG, \bphi_m \rangle \in \LL^p(0,T)$ for all $p \in \big[1, \frac{1}{1-\alpha}\big)$, in conjunction with the initial conditions
\[ \beta_m(0) = (\varrho \bg,\bphi_m),\qquad \dot{\beta}_m(0) = (\varrho \bh,\bphi_m),\qquad \mbox{$m=1,\dots,n$}.\]
The existence of a unique solution $\beta_m$ to this problem, with $\dot{\beta}_m  \in \mathrm{AC}([0,T])$ for each $m \in \{1,\dots,n\}$ is easily shown:
by letting $\gamma_m:= \dot{\beta}_m$, \eqref{eq:27} can be rewritten as a first-order system
for the two-component function $t \in [0,T] \mapsto (\beta_m(t),\gamma_m(t))^{\rm T} \in \mathbb{R}^2$, and then, because $\langle \bG, \bphi_m \rangle
\in \LL^1(0,T)$, integration of this system
over $[0,t]$, with $t \in (0,T]$ yields an integral equation
to which one can apply Banach's fixed point theorem in the complete metric space $\CC([0,T]) \times \CC([0,T])$ to deduce the existence of a unique absolutely continuous solution $(\beta_m, \gamma_m)^{\rm T}$, defined
on a ``maximal'' interval $[0,t_*] \subset [0,T]$. If $t_*$ were strictly less than $T$,
then it would follow that $|\beta_m(t)| + |\gamma_m(t)| \rightarrow + \infty$ as $t \rightarrow t_*$; the \textit{a priori} bound \eqref{eq:29}, which we shall prove below, however rules out this possibility; therefore $t_* = T$. Thus we deduce the existence of a unique Galerkin approximation $t \in [0,T] \mapsto \bu_n(t) \in \mathcal{V}_n$, with $\dot{\bu}_n \in \mathrm{AC}([0,T]; \mathcal{V}_n)$.

By taking $\bv=\bphi_m$ in \eqref{eq:26}, multiplying the resulting equality with $\beta_m(t)$ and summing over $m=1,\dots,n$,
we deduce that
\begin{align*}
&\frac{\tau^\alpha}{2} \frac{\dd}{\dd t}\int_\Omega \varrho |\dot \bu_n(t,\bx)|^2 \dd \bx  + (1-\tau^\alpha)\, \int_\Omega \frac{\partial}{\partial t}(-\dot{e}_{\alpha,1}\ast_t \sqrt{\varrho} \dot{\bu}_n)(t,\bx)\cdot \sqrt{\varrho}\dot{\bu}_n(t,\bx) \dd \bx \\
&\qquad + \frac{1}{2} \frac{\dd}{\dd t} \int_\Omega 2\mu |\beps(\bu_n(t,\bx))|^2 + \lambda |\tr(\beps(\bu_n(t,\bx)))|^2  \dd \bx \, = \, \langle \bG(t,\cdot), \dot{\bu}_n(t,\cdot)\rangle.
\end{align*}
Hence, by integration over $t \in (0,T]$ and noting the initial conditions satisfied by $\bu_n$ we deduce that
\begin{align*}
&\frac{\tau^\alpha }{2} \int_\Omega \varrho |\dot \bu_n(t,\bx)|^2 \dd \bx  + (1-\tau^\alpha)\, \int_0^t \int_\Omega \frac{\partial}{\partial s}(-\dot{e}_{\alpha,1}\ast_s \sqrt{\varrho} \dot{\bu}_n)(s,\bx)\cdot \sqrt{\varrho} \dot{\bu}_n(s,\bx) \dd \bx \dd s \\
&\quad \quad + \frac{1}{2} \int_\Omega 2\mu |\beps(\bu_n(t,\bx))|^2 + \lambda |\tr(\beps(\bu_n(t,\bx)))|^2  \dd \bx \\
&= \, \int_0^t  \langle \bG(s,\cdot), \dot{\bu}_n(s,\cdot)\rangle \dd s \\
&\quad \quad + \frac{\tau^\alpha }{2} \int_\Omega \varrho |\dot{\bu}_n(0,\bx)|^2 \dd \bx + \frac{1}{2} \int_\Omega 2\mu |\beps(\bu_n(0,\bx))|^2 + \lambda |\tr(\beps(\bu_n(0,\bx)))|^2  \dd \bx
\\
&= \, \int_0^t \langle \bG(s,\cdot), \dot{\bu}_n(s,\cdot)\rangle  \dd s + \frac{\tau^\alpha }{2} \|\dot{\bu}_n(0,\cdot)\|^2_{\LL^2_\varrho(\Omega)}
+ \frac{1}{2} a(\bu_n(0,\cdot), \bu_n(0,\cdot))\\
&\leq \, \int_0^t \langle \bG(s,\cdot), \dot{\bu}_n(s,\cdot)\rangle  \dd s + \frac{\tau^\alpha }{2} \|\bh\|^2_{\LL^2_\varrho(\Omega)}
+ \frac{1}{2} a(\bg, \bg)\\
&= \int_0^t \langle \bG(s,\cdot), \dot{\bu}_n(s,\cdot)\rangle  \dd s + \frac{\tau^\alpha }{2} \int_\Omega \varrho |\bh|^2 \dd \bx + \frac{1}{2} \int_\Omega 2\mu |\beps(\bg(\bx))|^2 + \lambda |\tr(\beps(\bg(\bx)))|^2  \dd \bx.
\end{align*}
Therefore,
\begin{align}\label{eq:28}
\begin{aligned}
&\frac{\tau^\alpha }{2} \int_\Omega \varrho |\dot \bu_n(t,\bx)|^2 \dd \bx  + (1-\tau^\alpha)\, \int_0^t \int_\Omega \frac{\partial}{\partial s}(-\dot{e}_{\alpha,1}\ast_s \sqrt{\varrho}\dot{\bu}_n)(s,\bx)\cdot \sqrt{\varrho}\dot{\bu}_n(s,\bx) \dd \bx \dd s \\
&\quad \quad + \frac{1}{2} \int_\Omega 2\mu |\beps(\bu_n(t,\bx))|^2 + \lambda |\tr(\beps(\bu_n(t,\bx)))|^2  \dd \bx \\
& \leq \int_0^t \langle \bG(s,\cdot), \dot{\bu}_n(s,\cdot)\rangle  \dd s + \frac{\tau^\alpha }{2} \int_\Omega \varrho |\bh|^2 \dd \bx + \frac{1}{2} \int_\Omega 2\mu |\beps(\bg(\bx))|^2 + \lambda |\tr(\beps(\bg(\bx)))|^2  \dd \bx.
\end{aligned}
\end{align}
We can now repeat the procedure (this time rigorously, as $\bu_n$ possesses the necessary regularity properties) leading from \eqref{eq:9} to the energy inequality \eqref{eq:18}, with $\bu$ replaced by $\bu_n$ throughout, resulting in the uniform bound
\begin{align}\label{eq:29}
\begin{aligned}
&\frac{\tau^\alpha}{2} \|\dot \bu_n(t)\|^2_{\LL^2_\varrho (\Omega)} + \frac{1}{2} \|\beps(\bu_n(t))\|^2_{\LL^2_\mu(\Omega)} + \frac{1}{2} \|\tr(\beps(\bu_n(t)))\|^2_{\LL^2_\lambda(\Omega)}\\
&\quad + (1-\tau^\alpha)\, \int_0^t \int_\Omega \frac{\partial}{\partial s}(-\dot{e}_{\alpha,1}\ast_s \sqrt{\varrho}\dot{\bu}_n)(s)\cdot \sqrt{\varrho}\dot{\bu}_n(s) \dd \bx \dd s \leq A(t) \, \mathrm{exp}(t+1 - e_{\alpha,1}(t)),
\end{aligned}
\end{align}
for all $t \in (0,T]$, with $A(t)$ again defined by the expression \eqref{eq:17}.

We are now ready to pass to the limit $n \rightarrow \infty$. To this end,
we fix an integer $N$ and choose a function  $\bv \in \CC^2_0([0,T);[\HH^1_0(\Omega)]^3)$ of the form
\begin{align}\label{eq:30}
\bv(t,\bx) := \sum_{k=1}^N \alpha_k(t) \bphi_k(\bx),
\end{align}
where $\alpha_k \in \CC^2_0([0,T))$ for $k=1,\dots,N$, i.e., $\alpha_k \in \CC^2([0,T])$ and
has compact support in the half-open interval $[0,T)$.
We then choose $n \geq N$ in \eqref{eq:26},
take $\bv = \bphi_k$ as test function in \eqref{eq:26} for $k \in \{1,\dots,N\}$, multiply the
resulting equality with $\alpha_k$, sum through $k=1,\dots,N$, and perform partial integrations in the first and the second term on the left-hand side to deduce that
\begin{align*}
&\tau^\alpha (\varrho \bu_n (0,\cdot), \dot \bv(0,\cdot)) - \tau^\alpha (\varrho \dot \bu_n (0,\cdot), \bv(0,\cdot))
+ \tau^\alpha \int_0^T (\varrho \bu_n (s,\cdot), \ddot{\bv}(s,\cdot) ) \dd s \\
&- (1-\tau^\alpha) ((-\dot{e}_{\alpha,1}\ast_t \varrho \dot{\bu}_n)(0,\cdot), \bv(0,\cdot)) - (1-\tau^\alpha) \int_0^T ((-\dot{e}_{\alpha,1}\ast_s \varrho \dot{\bu}_n)(s,\cdot), \dot{\bv}(s,\cdot)) \dd s \\
&+ \int_0^T \big ( 2\mu \beps(\bu_n(s,\cdot)) + \lambda \tr(\beps(\bu_n(s,\cdot))) \bI
\,,  \beps(\bv(s,\cdot)) \big) \dd s = \int_0^T \langle \bG(s,\cdot), \bv(s,\cdot) \rangle \dd s
\end{align*}
for all $\bv$ as in \eqref{eq:30} with $N$ fixed, and with any $n \geq N$.

Thus, because $(\varrho \bu_n (0,\cdot), \dot \bv(0,\cdot))=(\varrho \bg, \dot \bv(0,\cdot))$ and $(\varrho \dot \bu_n (0,\cdot), \bv(0,\cdot)) = (\varrho \bh, \bv(0,\cdot))$ for all $\bv \in \mathcal{V}_n$, and therefore (since $n \geq N$) also for all $\bv$ of the form \eqref{eq:30}, and as $ ((-\dot{e}_{\alpha,1}\ast_t \varrho \dot{\bu}_n)(0,\cdot) =0$, we have that
\begin{align}\label{eq:31}
\begin{aligned}
 \tau^\alpha &\int_0^T (\varrho \bu_n (s,\cdot), \ddot{\bv}(s,\cdot)) \dd s - (1-\tau^\alpha) \int_0^T ((-\dot{e}_{\alpha,1}\ast_s \varrho \dot{\bu}_n)(s,\cdot), \dot{\bv}(s,\cdot)) \dd s \\
&\quad + \int_0^T \big ( 2\mu \beps(\bu_n(s,\cdot)) + \lambda \tr(\beps(\bu_n(s,\cdot))) \bI
\,,  \beps(\bv(s,\cdot)) \big) \dd s \\
& = - \tau^\alpha  (\varrho \bg, \dot{\bv}(0,\cdot)) + \tau^\alpha  (\varrho \bh, \bv(0,\cdot)) + \int_0^T \langle \bG(s,\cdot), \bv(s,\cdot) \rangle \dd s.
\end{aligned}
\end{align}

As $0 \leq A(t) \leq A(T)=:A$ and $\mathrm{exp}(t+1 - e_{\alpha,1}(t)) \leq \mathrm{exp}(T+1)$, it follows from the energy estimate \eqref{eq:29} and Lemma \ref{le:1} that
\begin{itemize}
\item $(\bu_n)_{n \in \mathbb{N}}$ is a bounded sequence in $\LL^\infty(0,T;[\HH^1_0(\Omega)]^3)$;
\item $(\dot \bu_n)_{n \in \mathbb{N}}$ is a bounded sequence in $\LL^\infty(0,T;[\LL^2_\varrho(\Omega)]^3)\simeq \LL^\infty(0,T;[\LL^2(\Omega)]^3)$;
\item $((-\dot e_{\alpha,1})^{\frac{1}{2}}\,\dot \bu_n)_{n \in \mathbb{N}}$ is a bounded sequence in $\LL^2(0,T;[\LL^2_\varrho (\Omega)]^3) \simeq \LL^2(0,T;[\LL^2(\Omega)]^3)$.
\end{itemize}
%
%
Thus, by the Banach--Alaoglu theorem there exists a subsequence $(\bu_{n_\ell})_{\ell=1}^\infty$ such that
\begin{align}\label{eq:32}
 \left\{ \begin{array}{rll} \bu_{n_\ell} &\rightharpoonup \;\;\bu & \qquad \mbox{weakly$^\ast$ in $\LL^\infty(0,T;[\HH^1_0(\Omega)]^3)$}, \\
  \dot{\bu}_{n_\ell} &\rightharpoonup \;\;\dot{\bu} & \qquad \mbox{weakly$^\ast$ in $\LL^\infty(0,T;[\LL^2(\Omega)]^3)$},\\
  (-\dot e_{\alpha,1})^{\frac{1}{2}}\,\dot{\bu}_{n_\ell} &\rightharpoonup \;\;(-\dot e_{\alpha,1})^{\frac{1}{2}}\,\dot{\bu} & \qquad \mbox{weakly in $\LL^2(0,T;[\LL^2(\Omega)]^3)$}.
         \end{array}
 \right.
\end{align}
Furthermore, because for any $s \in (\frac{1}{2},1)$ the Sobolev space $[\HH^1_0(\Omega)]^3$ is compactly embedded into the fractional-order Sobolev space $[\HH^s_0(\Omega)]^3$, which is, in turn, continuously embedded into $[\LL^2(\Omega)]^3$, it follows from the Aubin--Lions--Simon lemma (cf. \cite{BS}) and the first two bullet points above that
\begin{align}\label{eq:33}
 \begin{array}{rll}
          \bu_{n_\ell} &\rightarrow \;\;\bu & \qquad \mbox{strongly in $\CC([0,T];[\HH^s_0(\Omega)]^3)$},\qquad s \in (\frac{1}{2},1),
         \end{array}
\end{align}
and therefore also
\begin{align}\label{eq:34}
 \begin{array}{rll}
          \bu_{n_\ell} &\rightarrow \;\;\bu & \qquad \mbox{strongly in $\CC([0,T];[\LL^2(\Omega)]^3)$}.
         \end{array}
\end{align}

We take $n=n_\ell$ in \eqref{eq:31} and pass to the limit $\ell \rightarrow \infty$ with
$\bv$ fixed. It then follows that
\begin{align}\label{eq:35}
\begin{aligned}
\tau^\alpha &\int_0^T (\varrho\bu (s,\cdot), \ddot{\bv}(s,\cdot)) \dd s - (1-\tau^\alpha) \int_0^T ((-\dot{e}_{\alpha,1}\ast_s \varrho \dot{\bu})(s,\cdot), \dot{\bv}(s,\cdot)) \dd s \\
&\quad + \int_0^T \big ( 2\mu \beps(\bu(s,\cdot)) + \lambda \tr(\beps(\bu(s,\cdot))) \bI
\,,  \beps(\bv(s,\cdot)) \big) \dd s \\
& = -\tau^\alpha  (\varrho\bg, \dot{\bv}(0,\cdot)) + \tau^\alpha  (\varrho\bh, \bv(0,\cdot)) + \int_0^T \langle \bG(s,\cdot), \bv(s,\cdot) \rangle \dd s
\end{aligned}
\end{align}
for all $\bv$ as in \eqref{eq:30} above, with $N$ fixed. This equality however holds
for all functions $\bv \in \WW^{2,1}(0,T;[\LL^2(\Omega)]^3) \cap \LL^1(0,T;[\HH^1_0(\Omega))]^3)$
such that $\bv(T,\cdot)=0$ and $\dot{\bv}(T,\cdot)=0$, as the set of all functions of the form \eqref{eq:30} is dense in this function space.

We note here that the passage to the limit in the second term on the left-hand side of \eqref{eq:31}, resulting in the second term on the left-hand side of \eqref{eq:35}  proceeds as follows: by Fubini's theorem to interchange the spatial integral with the integral with respect to $s$, and then by interchanging the order of integration in $s$ and $t$,
we have that
\begin{align*}
 \int_0^T ((-\dot{e}_{\alpha,1}\ast_s \varrho\dot{\bu}_{n_\ell})(s,\cdot), \dot{\bv}(s,\cdot)) \dd s &= - \int_0^T \left(\int_0^t \dot{e}_{\alpha,1}(t-s) \varrho\dot{\bu}_{n_\ell}(s,\cdot) \dd s , \dot{\bv}(t,\cdot)\right)\! \dd t\\
& = - \int_0^T \int_0^t \dot{e}_{\alpha,1}(t-s) \left(\varrho\dot{\bu}_{n_\ell}(s,\cdot), \dot{\bv}(t,\cdot)\right)\! \dd s \dd t \\
& = - \int_0^T \int_s^T \dot{e}_{\alpha,1}(t-s) \left(\varrho\dot{\bu}_{n_\ell}(s,\cdot), \dot{\bv}(t,\cdot)\right)\! \dd t \dd s \\
& = - \int_0^T \left(\varrho\dot{\bu}_{n_\ell}(s,\cdot), \int_s^T \dot{e}_{\alpha,1}(t-s) \dot{\bv}(t,\cdot)\dd t \right) \! \dd s.
\end{align*}
Then, because
$s \in [0,T] \mapsto \int_s^T \dot{e}_{\alpha,1}(t-s) \dot{\bv}(t,\cdot)\dd t \in
 \LL^1(0,T;[\LL^2(\Omega)]^3),
$
noting \eqref{eq:32}$_2$ yields
\begin{align*}
\lim_{\ell \rightarrow \infty} \int_0^T ((-\dot{e}_{\alpha,1}\ast_s \varrho\dot{\bu}_{n_\ell})(s,\cdot), \dot{\bv}(s,\cdot)) \dd s &= - \lim_{\ell \rightarrow \infty} \int_0^T \left(\varrho\dot{\bu}_{n_\ell}(s,\cdot), \int_s^T \dot{e}_{\alpha,1}(t-s) \dot{\bv}(t,\cdot)\dd t \right) \! \dd s\\
&= -\int_0^T \left(\varrho\dot{\bu}(s,\cdot), \int_s^T \dot{e}_{\alpha,1}(t-s) \dot{\bv}(t,\cdot)\dd t \right)\!  \dd s\\
& = \int_0^T ((-\dot{e}_{\alpha,1}\ast_t \varrho \dot{\bu}(s,\cdot), \dot{\bv}(s,\cdot))  \dd s,
\end{align*}
as has been asserted above. The passages to the limits in the first and third term on the left-hand side of \eqref{eq:31} are immediate, by using \eqref{eq:32}$_2$ and \eqref{eq:32}$_1$, respectively.

We have thereby shown the existence of a function $\bu \in \LL^\infty(0,T;[\HH^1_0(\Omega)]^3)$
such that $\dot{\bu} \in \LL^\infty(0,T;[\LL^2(\Omega)]^3)$, satisfying \eqref{eq:35} for all
$\bv \in \WW^{2,1}(0,T;[\LL^2(\Omega)]^3) \cap \LL^1(0,T;[\HH^1_0(\Omega))]^3)$
such that $\bv(T,\cdot)=0$ and $\dot{\bv}(T,\cdot)=0$;
the proof of the existence of a weak solution is therefore almost complete. It remains
to show that $\bu \in \CC_w([0,T];[\HH^1_0(\Omega)]^3)$
and $\dot{\bu} \in \CC_w([0,T];[\LL^2(\Omega)]^3)$.

We begin by recalling that, for any pair of Hilbert spaces $\mathcal{H}$ and $\mathcal{V}$ such that $\mathcal{V}$ is continuously and densely embedded into $\mathcal{H}$,
if $v \in \LL^\infty(0,T;\mathcal{V})$ and $\dot{v} \in \LL^1(0,T;\mathcal{H})$ (whereby $v \in \WW^{1,1}(0,T;\mathcal{H})
\subset \CC([0,T];\mathcal{H}) \subset \CC_w([0,T];\mathcal{H})$),
then $v \in \CC_w([0,T];\mathcal{V})$ (cf. eq. (8.49) in Lemma 8.1, Ch. 3 of \cite{LM}). Therefore, because
\[\bu \in \LL^\infty(0,T;[\HH^1_0(\Omega)]^3)\quad\mbox{and}\quad
%
\dot{\bu} \in \LL^\infty(0,T;[\LL^2(\Omega)]^3) \subset \LL^1(0,T;[\LL^2(\Omega)]^3),
\]
it follows, with $\mathcal{V}=[\HH^1_0(\Omega)]^3$ and $\mathcal{H} = [\LL^2(\Omega)]^3$, that $\bu \in \CC_w([0,T];[\HH^1_0(\Omega)]^3)$.

Next, we will show that $\dot{\bu} \in \CC_w([0,T];[\LL^2(\Omega)]^3)$. It follows from \eqref{eq:26} that
$\bu_n(t) \in \mathcal{V}_n$, for $t \in [0,T]$, satisfies:
\begin{align}\label{eq:36}
\hspace{-3mm} \left( \frac{\partial}{\partial t}\big(\tau^\alpha  \varrho \dot \bu_n + (1-\tau^\alpha)
(-\dot{e}_{\alpha,1}\ast_t \varrho \dot{\bu}_n)\big), \bv \right ) = - \big ( 2\mu \beps(\bu_n) + \lambda \tr(\beps(\bu_n)) \bI
\,, \beps(\bv) \big) + \langle \bG, \bv \rangle
\end{align}
for all $\bv \in \mathcal{V}_n$. We thus have from \eqref{eq:36} that,
for any $\bv \in [\HH^1_0(\Omega)]^3$,
\begin{align*}
\hspace{-3mm} \left( \frac{\partial}{\partial t}\big(\tau^\alpha \varrho \dot \bu_n + (1-\tau^\alpha)
(-\dot{e}_{\alpha,1}\ast_t \varrho \dot{\bu}_n)\big), \bv \right ) &= \left( \frac{\partial}{\partial t}\big(\varrho \tau^\alpha\dot \bu_n + \varrho (1-\tau^\alpha) (-\dot{e}_{\alpha,1}\ast_t \dot{\bu}_n)\big), P_n \bv \right )\\
& =  - \big ( 2\mu \beps(\bu_n) + \lambda \tr(\beps(\bu_n)) \bI
\,, \beps(P_n \bv) \big) + \langle \bG, P_n \bv \rangle.
\end{align*}
We note that by the energy estimate \eqref{eq:29} and because $\|P_n\|_{\mathcal{L}([\HH^1(\Omega)]^3, [\HH^1(\Omega)]^3)}$ is bounded by $(c_1/c_0)^{1/2}$, uniformly in $n$, there exists a positive constant $C$, independent of $n$, such that
\begin{alignat*}{2}
&\big ( 2\mu \beps(\bu_n) + \lambda \tr(\beps(\bu_n)) \bI
\,, \beps(P_n \bv) \big) =   \big (2\mu \beps(\bu_n)\,,  \beps(P_n \bv) \big) + \big(\lambda \tr(\beps(\bu_n)) \bI
\,, \beps(P_n \bv) \big)\\
&\qquad =   \big (2\mu \beps(\bu_n)\,,  \beps(P_n \bv) \big) + \big(\lambda \tr(\beps(\bu_n))
\,, \tr(\beps(P_n \bv)) \big)\\
&\qquad \leq \left(2\|\beps(\bu_n)\|^2_{\LL^2_\mu(\Omega)} + \|\tr(\beps(\bu_n))\|^2_{\LL^2_\lambda(\Omega)}\right)^{\frac{1}{2}}
\left(2\|\beps(P_n \bv)\|^2_{\LL^2_\mu(\Omega)} + \|\tr(\beps(P_n \bv))\|^2_{\LL^2_\lambda(\Omega)}\right)^{\frac{1}{2}}\\
&\qquad \leq C \left(2\|\beps(P_n \bv)\|^2_{\LL^2_\mu(\Omega)} + \|\tr(\beps(P_n \bv))\|^2_{\LL^2_\lambda(\Omega)}\right)^{\frac{1}{2}} \leq C  \|\bv \|_{\HH^1(\Omega)}
\qquad
\forall\, \bv \in [\HH^1_0(\Omega)]^3.
\end{alignat*}
Also,
\begin{align*}
\langle \bG , P_n \bv \rangle & = (\tau^\alpha-1) \dot{e}_{\alpha,1} \left(\varrho \bh\,, P_n \bv \right)  \,-\, e_{\alpha,1} \left(\tau^\alpha\bs-
2\mu \beps(\bg) - \lambda \tr(\beps(\bg)) \bI\,, \beps(P_n \bv) \right)\\
& \qquad +  \tau^\alpha \left(\bF , P_n \bv\right) + (\tau^\alpha-1) \dot{e}_{\alpha,1} \ast_t \left(\bF , P_n \bv\right) \qquad
\forall\, \bv \in [\HH^1_0(\Omega)]^3,
\end{align*}
and therefore, because $\|P_n\|_{\mathcal{L}([\LL^2(\Omega)]^3, [\LL^2(\Omega)]^3)}\leq 1$ and $\|P_n\|_{\mathcal{L}([\HH^1(\Omega)]^3, [\HH^1(\Omega)]^3)}$ is bounded by $(c_1/c_0)^{1/2}$, uniformly in $n$, there exists a positive constant $C$, independent of $n$, such that
\begin{align*}
|\langle \bG , P_n \bv \rangle| & \leq (1-\tau^\alpha) (-\dot{e}_{\alpha,1}) \|\bh\|_{\LL^2_\varrho(\Omega)} \|\bv\|_{\LL^2_\varrho(\Omega)}\\
&\qquad + \, C e_{\alpha,1} \|\tau^\alpha\bs-
2\mu \beps(\bg) - \lambda \tr(\beps(\bg)) \bI\|_{\LL^2(\Omega)} \|\bv\|_{\HH^1(\Omega)}\\
& \qquad +  \tau^\alpha \|\bF\|_{\LL^2_{1/\varrho}(\Omega)} \|\bv\|_{\LL^2_{\varrho}(\Omega)}\\
& \qquad + (1-\tau^\alpha) \big((-\dot{e}_{\alpha,1}) \ast_t \|\bF\|_{\LL^2_{1/\varrho}(\Omega)}\big) \|\bv\|_{\LL^2_{\varrho}(\Omega)}
\qquad
\forall\, \bv \in [\HH^1_0(\Omega)]^3.
\end{align*}
Thus we deduce, with $\varrho_1:= \|\varrho\|_{\LL^\infty(\Omega)}$,
\begin{align*}
&\left\|\frac{\partial}{\partial t}\big(\tau^\alpha \varrho \dot \bu_n + (1-\tau^\alpha)
(-\dot{e}_{\alpha,1}\ast_t \varrho \dot{\bu}_n)\big)\right\|_{\HH^{-1}(\Omega)} \\
&:= \sup_{\bv \in [\HH^1_0(\Omega)]^3} \frac{\left\langle \frac{\partial}{\partial t}\big(\varrho \tau^\alpha\dot \bu_n + \varrho (1-\tau^\alpha)
(-\dot{e}_{\alpha,1}\ast_t \dot{\bu}_n)\big) , \bv \right\rangle}{\|\bv\|_{\HH^1_0(\Omega)}}\\
& \leq C  + \varrho_1 (1-\tau^\alpha) (-\dot{e}_{\alpha,1}) \|\bh\|_{\LL^2(\Omega)} \\
&\qquad +  C e_{\alpha,1}\|\tau^\alpha\bs-
2\mu \beps(\bg) - \lambda \tr(\beps(\bg)) \bI\|_{\LL^2(\Omega)} \\
&\qquad +  \tau^\alpha \sqrt{\frac{\varrho_1}{\varrho_0}} \|\bF\|_{\LL^2(\Omega)} + (1-\tau^\alpha) \sqrt{\frac{\varrho_1}{\varrho_0}} \big((-\dot{e}_{\alpha,1}) \ast_t \|\bF\|_{\LL^2(\Omega)}\big),
\end{align*}
which then implies, because of \eqref{eq:23}, for any $p \in [1,2]$ satisfying $p<\frac{1}{1-\alpha}$, that
\begin{align*}
& \left\|\frac{\partial}{\partial t}\big(\tau^\alpha\varrho \dot \bu_n + \varrho (1-\tau^\alpha)
(-\dot{e}_{\alpha,1}\ast_t \varrho \dot{\bu}_n)\big)\right\|_{\LL^p(0,T;\HH^{-1}(\Omega))}  \\
&\leq  C T^{\frac{1}{p}} + \varrho_1 (1-\tau^\alpha) \|-\dot{e}_{\alpha,1}\|_{\LL^p(0,T)} \|\bh\|_{\LL^2(\Omega)}
\end{align*}
\begin{align*}
&\qquad + C \|e_{\alpha,1}\|_{\LL^p(0,T)} \|\tau^\alpha\bs-
2\mu \beps(\bg) - \lambda \tr(\beps(\bg)) \bI\|_{\LL^2(\Omega)} \\
&\qquad +  \tau^\alpha \sqrt{\frac{\varrho_1}{\varrho_0}} \|\bF\|_{\LL^p(0,T;\LL^2(\Omega))} + (1-\tau^\alpha) \sqrt{\frac{\varrho_1}{\varrho_0}} \|-\dot{e}_{\alpha,1}\|_{\LL^1(0,T)} \|\bF\|_{\LL^p(0,T;\LL^2(\Omega))}.
\end{align*}
Hence, for any $p \in [1,2]$ such that $p<\frac{1}{1-\alpha}$, we have that
\begin{align*}
\left\|\frac{\partial}{\partial t}\big(\tau^\alpha \varrho \dot \bu_n + \varrho (1-\tau^\alpha)
(-\dot{e}_{\alpha,1}\ast_t \varrho \dot{\bu}_n)\big) \right\|_{\LL^p(0,T;\HH^{-1}(\Omega))} \leq C,
\end{align*}
where $C$ is a positive constant, independent of $n$. Consequently, by the Banach--Alaoglu theorem, there exists a subsequence $(\bu_{n_\ell})_{\ell=1}^\infty$ such that
\[\frac{\partial}{\partial t}\big(\tau^\alpha \varrho \dot \bu_{n_\ell} + (1-\tau^\alpha)
(-\dot{e}_{\alpha,1}\ast_t \varrho\dot{\bu}_{n_\ell})\big) \rightharpoonup \frac{\partial}{\partial t}\big(\tau^\alpha\varrho \dot \bu + (1-\tau^\alpha)
(-\dot{e}_{\alpha,1}\ast_t \varrho \dot{\bu})\big), \]
weakly in $\LL^p(0,T;[\HH^{-1}(\Omega)]^3)$ for any $p \in [1,2]$ such that $p<\frac{1}{1-\alpha}$. As
\[ \tau^\alpha \varrho\dot \bu + (1-\tau^\alpha) (-\dot{e}_{\alpha,1}\ast_t \varrho\dot{\bu}) \in \LL^\infty(0,T;[\LL^2(\Omega)]^3)\]
and
\[ \frac{\partial}{\partial t}\big(\tau^\alpha\varrho \dot \bu + (1-\tau^\alpha)
(-\dot{e}_{\alpha,1}\ast_t \varrho\dot{\bu})\big) \in \LL^1(0,T;[\HH^{-1}(\Omega)]^3),\]
it once again follows, thanks to the continuous embedding of $[\LL^2(\Omega)]^3$ into $[\HH^{-1}(\Omega)]^3$,  that
\begin{align}\label{eq:37}
\tau^\alpha \varrho \dot \bu + (1-\tau^\alpha) (-\dot{e}_{\alpha,1}\ast_t \varrho \dot{\bu}) \in \CC_w([0,T];[\LL^2(\Omega)]^3).
\end{align}
However, as $\varrho \dot{\bu} \in \LL^\infty(0,T;[\LL^2(\Omega)]^3)$, we have that $t \in [0,T] \mapsto (\varrho\dot{\bu}(t),\bw)$
belongs to $\LL^\infty(0,T)$ for each $\bw \in [\LL^2(\Omega)]^3$, and therefore, thanks to the smoothing property of the convolution, the function
$t \in [0,T] \mapsto -\dot{e}_{\alpha,1}(t) \ast_t (\varrho \dot{\bu}(t),\bw)$ belongs to $\CC([0,T])$. Consequently,
\[ t \in [0,T] \mapsto (-\dot{e}_{\alpha,1}(t) \ast_t \varrho \dot{\bu}(t),\bw) \in \CC([0,T])\qquad \forall\, \bw \in [\LL^2(\Omega)]^3,\]
meaning that $(1-\tau^\alpha) (-\dot{e}_{\alpha,1}) \ast_t \varrho  \dot{\bu} \in \CC_w([0,T];[\LL^2(\Omega)]^3)$, and therefore
by \eqref{eq:37}, also $\varrho\dot{\bu} \in \CC_w([0,T];[\LL^2(\Omega)]^3)$. Because $\varrho_0 \leq \varrho(\bx) \leq \varrho_1$ a.e. on $\Omega$, it then follows that
\[\dot{\bu} \in \CC_w([0,T];[\LL^2(\Omega)]^3).\]
That completes the proof of the existence of a weak solution.

\smallskip

STEP 2: \textit{Proof of the energy inequality.} Next we prove that weak solutions whose existence we have thus proved satisfy the energy inequality in the statement of the
theorem. Our starting point is \eqref{eq:29}.
By Lemma \ref{le:1}, we have that
\begin{align*}
&\int_0^t \int_\Omega \frac{\partial}{\partial s}(-\dot{e}_{\alpha,1}\ast_s \sqrt{\varrho}\dot{\bu}_n)(s)\cdot \sqrt{\varrho} \dot{\bu}_n(s) \dd \bx \dd s
\\
&\qquad \geq  \frac{1}{2} (-\dot{e}_{\alpha,1}(\cdot) \ast_t \|\sqrt{\varrho}\dot{\bu}_n(\cdot)\|^2_{\LL^2(\Omega)})(t)
+ \frac{1}{2}\int_0^t -\dot{e}_{\alpha,1}(s)\|\sqrt{\varrho}\dot{\bu}_n(s)\|^2_{\LL^2(\Omega)} \dd s
\qquad \mbox{for all $t \in (0,T]$},
\end{align*}
and each of the two terms on the right-hand side is nonnegative. By omitting the first term
from the right-hand side of this equality, and substituting the resulting inequality into \eqref{eq:29} we have that
\begin{align}\label{eq:39}
\begin{aligned}
&\frac{\tau^\alpha}{2} \|\dot \bu_n(t)\|^2_{\LL^2_\varrho(\Omega)} + \frac{1}{2} \|\beps(\bu_n(t))\|^2_{\LL^2_\mu(\Omega)} + \frac{1}{2} \|\tr(\beps(\bu_n(t)))\|^2_{\LL^2_\lambda(\Omega)}\\
&\quad + \frac{1-\tau^\alpha}{2} \int_0^t -\dot{e}_{\alpha,1}(s)\|\dot{\bu}_n(s)\|^2_{\LL^2_\varrho(\Omega)} \dd s
\leq A(t) \, \mathrm{exp}(t+1 - e_{\alpha,1}(t))\quad \mbox{for all $t \in (0,T]$}.
\end{aligned}
\end{align}

As $\dot{\bu}_{n_\ell} \rightharpoonup \dot{\bu}$ weakly$^\ast$ in $\LL^\infty(0,T;[\LL^2_\varrho(\Omega)]^3)$, the weak lower-semicontinuity of the norm function and \eqref{eq:39} imply that
\begin{align}\label{eq:40}
\|\dot \bu(s)\|^2_{\LL^2_\varrho(\Omega)} \leq
\|\dot \bu\|^2_{\LL^\infty(0,t;\LL^2_\varrho(\Omega))} \leq \liminf_{\ell \rightarrow \infty} \|\dot \bu_{n_\ell}\|^2_{\LL^\infty(0,t;\LL^2_\varrho(\Omega))} \leq A(t) \, \mathrm{exp}(t+1 - e_{\alpha,1}(t))
\end{align}
for all $t \in (0,T]$ and a.e. $s \in (0,t]$. Similarly, because $\bu_{n_\ell} \rightharpoonup \bu$ weakly$^\ast$ in $\LL^\infty(0,T;[\HH^1_0(\Omega)]^3)$,
\begin{align}\label{eq:41}
\begin{aligned}
&\bigg[\frac{1}{2} \|\beps(\bu(s))\|^2_{\LL^2_\mu(\Omega)} + \frac{1}{2} \|\tr(\beps(\bu(s)))\|^2_{\LL^2_\lambda(\Omega)}\bigg]\\
&\qquad \leq \mbox{ess.sup}_{s \in (0,t]}\bigg[\frac{1}{2} \|\beps(\bu(s))\|^2_{\LL^2_\mu(\Omega)} + \frac{1}{2} \|\tr(\beps(\bu(s)))\|^2_{\LL^2_\lambda(\Omega)}\bigg]
\\
&\qquad \leq \liminf_{\ell \rightarrow \infty} \left\{
\mbox{ess.sup}_{s \in (0,t]}\left[\frac{\mu}{2} \|\beps(\bu_{n_\ell}(s))\|^2_{\LL^2(\Omega)} + \frac{1}{2} \|\tr(\beps(\bu_{n_\ell}(s)))\|^2_{\LL^2_\lambda(\Omega)}\right]\right\}\\
&\qquad \leq A(t) \, \mathrm{exp}(t+1 - e_{\alpha,1}(t))
\end{aligned}
\end{align}
for all $t \in (0,T]$ and a.e. $s \in (0,t]$. Finally, because
$(-e_{\alpha,1})^{\frac{1}{2}} \dot{\bu}_{n_\ell} \rightharpoonup (-e_{\alpha,1})^{\frac{1}{2}} \dot{\bu}$
weakly in {\color{black} the function space} $\LL^2(0,T;[\LL^2_\varrho(\Omega)]^3)$, we have that
\begin{align}\label{eq:42}
\int_0^t -\dot{e}_{\alpha,1}(s)\|\dot{\bu}(s)\|^2_{\LL^2_\varrho(\Omega)} \dd s \leq \liminf_{\ell \rightarrow \infty}  \int_0^t -\dot{e}_{\alpha,1}(s)\|\dot{\bu}_{n_\ell}(s)\|^2_{\LL^2_\varrho(\Omega)} \dd s \leq A(t) \, \mathrm{exp}(t+1 - e_{\alpha,1}(t))
\end{align}
for all $t \in (0,T]$. Summing \eqref{eq:40}--\eqref{eq:42} we deduce the asserted energy inequality \textcolor{black}{\eqref{EE}}.

\smallskip

{\color{black}
STEP 3: \textit{Attainment of the initial conditions for $\bu$ and $\dot{\bu}$.} Next, we shall prove that the initial condition $\bu(0,\cdot) = \bg(\cdot)$ is satisfied in the sense of continuous functions from $[0,T]$ into $[\LL^2(\Omega)]^3$. To this end, we note that
\begin{align*}
\|\bu(0,\cdot) - \bu_{n_\ell}(0,\cdot)\|_{\LL^2(\Omega)} \leq \|\bu - \bu_{n_\ell}\|_{\CC([0,T];\LL^2(\Omega))} \rightarrow 0 \qquad \mbox{as $\ell \rightarrow \infty$},
\end{align*}
thanks to \eqref{eq:34}. Since $\bu_{n_\ell}(0,\cdot) = P_{n_\ell}\bg(\cdot) \rightarrow \bg(\cdot)$ strongly in $[\LL^2_\varrho(\Omega)]^3\simeq [\LL^2(\Omega)]^3 $ as $\ell \rightarrow \infty$, we finally deduce by the triangle inequality that $\bu(0,\cdot) - \bg(\cdot) = 0$. Therefore, $\bu(0,\cdot) = \bg(\cdot)$, with $\bu \in \CC([0,T];[\LL^2(\Omega)]^3)$.

To show that the initial condition, $\dot{\bu}(0,\cdot) = \bh(\cdot)$ is satisfied we note that, thanks to \eqref{eq:20}$_1$ and \eqref{eq:20}$_2$,
we have $\bu \in \WW^{1,\infty}(0,T;[\LL^2(\Omega)]^3) = \CC^{0,1}([0,T];[\LL^2(\Omega)]^3)$, so we
can perform partial integration with respect to $t$ in the first term on the left-hand side of \eqref{eq:21}, resulting in
\begin{align*}
 -\tau^\alpha  (\varrho\bu(0,\cdot), \dot{\bv}(0,\cdot))&- \tau^\alpha \int_0^T (\varrho\dot\bu (s,\cdot), \dot{\bv}(s,\cdot)) \dd s
- (1-\tau^\alpha) \int_0^T ((-\dot{e}_{\alpha,1}\ast_s \varrho\dot{\bu})(s,\cdot), \dot{\bv}(s,\cdot)) \dd s \\
& + \int_0^T \big ( 2\mu \beps(\bu(s,\cdot)) + \lambda \tr(\beps(\bu(s,\cdot))) \bI
\,,  \beps(\bv(s,\cdot)) \big) \dd s \\
&\hspace{-5mm} = - \tau^\alpha  (\varrho\bg, \dot{\bv}(0,\cdot)) + \tau^\alpha  (\varrho\bh, \bv(0,\cdot)) + \int_0^T \langle \bG(s,\cdot), \bv(s,\cdot) \rangle \dd s
\end{align*}
for all $\bv \in \WW^{2,1}(0,T;[\LL^2(\Omega)]^3) \cap \LL^1(0,T;[\HH^1_0(\Omega)]^3)$
with $\bv(T,\cdot)=0$ and $\dot{\bv}(T,\cdot)=0$. As $\bu(0,\cdot) = \bg(\cdot)$, the first term on the left-hand side and the first term on the right-hand side cancel, whereby
\begin{align*}
&- \tau^\alpha \int_0^T (\varrho\dot\bu (s,\cdot), \dot{\bv}(s,\cdot)) \dd s
- (1-\tau^\alpha) \int_0^T ((-\dot{e}_{\alpha,1}\ast_s \varrho\dot{\bu})(s,\cdot), \dot{\bv}(s,\cdot)) \dd s \\
& + \int_0^T \big ( 2\mu \beps(\bu(s,\cdot)) + \lambda \tr(\beps(\bu(s,\cdot))) \bI
\,,  \beps(\bv(s,\cdot)) \big) \dd s
= \tau^\alpha  (\varrho\bh, \bv(0,\cdot)) + \int_0^T \langle \bG(s,\cdot), \bv(s,\cdot) \rangle \dd s
\end{align*}
for all $\bv \in \WW^{2,1}(0,T;[\LL^2(\Omega)]^3) \cap \LL^1(0,T;[\HH^1_0(\Omega)]^3)$
with $\bv(T,\cdot)=0$ and $\dot{\bv}(T,\cdot)=0$. As the set of all such $\bv$ is dense in the set of all
$\bv \in \WW^{1,1}(0,T;[\LL^2(\Omega)]^3) \cap \LL^1(0,T;[\HH^1_0(\Omega)]^3)$ with $\bv(T,\cdot)=0$, it follows that
\begin{align}\label{eq:4.23aa}
\begin{aligned}
&- \tau^\alpha \int_0^T (\varrho\dot\bu (s,\cdot), \dot{\bv}(s,\cdot)) \dd s
- (1-\tau^\alpha) \int_0^T ((-\dot{e}_{\alpha,1}\ast_s \varrho\dot{\bu})(s,\cdot), \dot{\bv}(s,\cdot)) \dd s\\
& + \int_0^T \big ( 2\mu \beps(\bu(s,\cdot)) + \lambda \tr(\beps(\bu(s,\cdot))) \bI
\,,  \beps(\bv(s,\cdot)) \big) \dd s
= \tau^\alpha  (\varrho\bh, \bv(0,\cdot)) + \int_0^T \langle \bG(s,\cdot), \bv(s,\cdot) \rangle \dd s
\end{aligned}
\end{align}
holds for all $\bv \in \WW^{1,1}(0,T;[\LL^2(\Omega)]^3) \cap \LL^1(0,T;[\HH^1_0(\Omega)]^3)$
with $\bv(T,\cdot)=0$.

We fix a $t_0 \in (0,T)$ and for $\varepsilon \in (0,T-t_0)$ we define
\[ \varphi_\varepsilon(t):= \left\{\begin{array}{cl} 1 & \mbox{for $0 \leq t \leq t_0$},\\
1 - \frac{1}{\varepsilon}(t - t_0) & \mbox{for $t_0 < t < t_0 + \varepsilon$},\\
0 & \mbox{for $t_0 + \varepsilon \leq t \leq T$.} \end{array} \right. \]
Clearly, $\varphi_\varepsilon \in \CC^{0,1}([0,T])$, the weak derivative of $\varphi_\varepsilon$ is $\varphi_\varepsilon' = - \frac{1}{\varepsilon} \chi_{(t_0,t_0 + \varepsilon)}$, and $\varphi_\varepsilon(T)=0$. Hence, for any $\bw \in [\HH^1(\Omega)]^3$, and taking $\bv = \varphi_\varepsilon \bw$ in \eqref{eq:4.23aa}, we have that
\begin{align}\label{eq:4.23bb}
&\tau^\alpha \,\frac{1}{\varepsilon}\int_{t_0}^{t_0+\varepsilon} (\varrho\dot\bu (s,\cdot), {\bw}(\cdot)) \dd s
+ (1-\tau^\alpha)\, \frac{1}{\varepsilon}\int_{t_0}^{t_0+\varepsilon} ((-\dot{e}_{\alpha,1}\ast_s \varrho\dot{\bu})(s,\cdot),{\bw}(\cdot)) \dd s \nonumber\\
& + \int_0^{t_0+\varepsilon} \big ( 2\mu \beps(\bu(s,\cdot)) + \lambda \tr(\beps(\bu(s,\cdot))) \bI
\,,  \varphi_\varepsilon(s)\, \beps(\bw(\cdot)) \big) \dd s\\
& \hspace{2in} = \tau^\alpha  (\varrho\bh, \bw) + \int_0^{t_0+\varepsilon} \langle \bG(s,\cdot), \varphi_\varepsilon(s) \bw(\cdot) \rangle \dd s \nonumber
\end{align}
for all $\bw \in [\HH^1_0(\Omega)]^3$. As $\varrho\dot{\bu} \in \CC_w([0,T];[\LL^2(\Omega)]^3)$ and $(-\dot{e}_{\alpha,1}) \ast_t \varrho  \dot{\bu} \in \CC([0,T];[\LL^2(\Omega)]^3)$ (cf. the end of STEP 1), we can pass to the limit $\varepsilon \rightarrow 0_+$ in \eqref{eq:4.23bb}, with $t_0 \in (0,T)$ fixed,  to deduce by applying Lebesgue's differentiation theorem to the first and the second integral on the left-hand side of \eqref{eq:4.23bb}, recalling the continuity of the integrands in those integrals as functions of the integration variable $s$, for $\bw \in [\HH^1_0(\Omega)]^3$ fixed, and using
the continuity of the integral with respect to its (upper) limit in the third integral on the left-hand side of
 \eqref{eq:4.23bb} and the second term on the right-hand side of \eqref{eq:4.23bb}, that
\begin{align}\label{eq:4.23dd}
\begin{aligned}
&\tau^\alpha \,(\varrho\dot\bu (t_0,\cdot), {\bw}(\cdot))
+ (1-\tau^\alpha)\, ((-\dot{e}_{\alpha,1}\ast_t \varrho\dot{\bu})(t_0,\cdot),{\bw}(\cdot))\\
& + \int_0^{t_0} \big ( 2\mu \beps(\bu(s,\cdot)) + \lambda \tr(\beps(\bu(s,\cdot))) \bI
\,,  \beps(\bw(\cdot)) \big) \dd s
= \tau^\alpha  (\varrho\bh, \bw)
+ \int_0^{t_0} \langle \bG(s,\cdot), \bw(\cdot) \rangle \dd s
\end{aligned}
\end{align}
for all $\bw \in [\HH^1_0(\Omega)]^3$ and all $t_0 \in (0,T)$. Next, with $\bw \in [\HH^1_0(\Omega)]^3$ fixed, we pass to the limit $t_0 \rightarrow 0_+$ in \eqref{eq:4.23dd}, noting that the third term on the left-hand side and the second term on the right-hand side both vanish in this limit thanks to the continuity of these integrals as functions of $t_0$, and that, for the same reason and by Fubini's theorem, also %
\[\lim_{t_0 \rightarrow 0_+} ((-\dot{e}_{\alpha,1}\ast_t \varrho\dot{\bu})(t_0,\cdot),{\bw}(\cdot)) = -\lim_{t_0 \rightarrow 0_+} \int_0^{t_0} \dot{e}_{\alpha,1}(t)\, (( \varrho\dot{\bu})(t - t_0,\cdot),{\bw}(\cdot)) \dd t = 0.\]
Consequently, upon passage to the limit $t_0 \rightarrow 0_+$, the equality \eqref{eq:4.23dd} collapses to
\[ \tau^\alpha \,(\varrho\dot\bu (0,\cdot), {\bw}(\cdot)) = \tau^\alpha  (\varrho\bh, \bw) \qquad \forall\, \bw \in [\HH^1_0(\Omega)]^3.\]
As $\tau \in (0,1]$ and $[\HH^1_0(\Omega)]^3$ is dense in $[\LL^2(\Omega)]^3$ it then follows that
\[ (\varrho\dot\bu (0,\cdot), {\bw}(\cdot)) = (\varrho\bh, \bw) \qquad \forall\, \bw \in [\LL^2(\Omega)]^3.\]
Because $\varrho \in \LL^\infty(\Omega)$ and $\varrho(x) \geq \rho_0>0$ a.e. in $\Omega$ (cf. \eqref{coeff-ass}), we
finally have that
\[ (\dot\bu (0,\cdot), {\bw}(\cdot)) = (\bh, \bw) \qquad \forall\, \bw \in [\LL^2(\Omega)]^3,\]
and therefore $\dot\bu (0,\cdot) = \bh(\cdot)$, as an equality in $\CC_w([0,T],[\LL^2(\Omega)]^3)$.
}

\smallskip

STEP 4. \textit{Uniqueness of the solution.}
Having shown, for
initial data $\bg$, $\bh$, $\bs$ and the source term $\bF$ satisfying \eqref{eq:19}, and for any
$\tau \in (0,1]$, $\varrho, \mu, \lambda \in \LL^\infty(\Omega)$, with $\varrho(\bx) \geq \varrho_0>0$, $\mu(\bx) \geq \mu_0>0$ and $\lambda(\bx) \geq 0$ for a.e. $\bx \in \Omega$, and $\alpha \in (0,1)$, the existence of a weak solution \eqref{eq:20}
to the problem \eqref{eq:2}, \eqref{eq:3}, \eqref{eq:8} satisfying the equality \eqref{eq:21} with \eqref{eq:22}
we now to turn to the proof of uniqueness of weak solutions.
{\color{black} Suppose that $\bu_1$ and $\bu_2$ are two weak solutions to \eqref{eq:2}, \eqref{eq:3}, \eqref{eq:8} subject to the same initial data and source term. Then, they both satisfy \eqref{eq:4.23aa}, and therefore, thanks to the linearity of the problem their difference $\bu:= \bu_1 - \bu_2 $ satisfies
\begin{align}\label{eq:46}
\begin{aligned}
- \tau^\alpha \int_0^T (\varrho\dot\bu (s,\cdot), \dot{\bv}(s,\cdot)) \dd s &- (1-\tau^\alpha) \int_0^T ((-\dot{e}_{\alpha,1}\ast_s \varrho\dot{\bu})(s,\cdot), \dot{\bv}(s,\cdot)) \dd s \\
& + \int_0^T \big ( 2\mu \beps(\bu(s,\cdot)) + \lambda \tr(\beps(\bu(s,\cdot))) \bI
\,,  \beps(\bv(s,\cdot)) \big) \dd s = 0
\end{aligned}
\end{align}
for all $\bv \in \WW^{1,1}(0,T;[\LL^2(\Omega)]^3) \cap \LL^1(0,T;[\HH^1_0(\Omega)]^3)$
with $\bv(T,\cdot)=0$. We fix a $t_0 \in (0,T)$,} and let
\[ \bv(t,\bx):= \left\{ \begin{array}{cl}
-\int_t^{t_0} \bu(s,\bx) \dd s & \mbox{for $0 < t \leq t_0$},\\
\mathbf{0} & \mbox{for $t_0 < t <T$.} \end{array} \right.\]
Clearly, $\bv \in  \WW^{1,\infty}(0,T;[\HH^1_0(\Omega)]^3)$
with $\bv(T,\cdot)=\mathbf{0}$, and hence the function $\bv$, thus defined, is an admissible test function. We therefore have from \eqref{eq:46} that
\begin{align}\label{eq:47}
\begin{aligned}
- \tau^\alpha \int_0^{t_0} (\varrho\dot\bu (s,\cdot), \bu(s,\cdot)) \dd s
&- (1-\tau^\alpha) \int_0^{t_0} ((-\dot{e}_{\alpha,1}\ast_s \varrho \dot{\bu})(s,\cdot),  \bu(s,\cdot)) \dd s \\
&+ \int_0^{t_0} \big ( 2\mu \beps(\dot\bv(s,\cdot)) + \lambda \tr(\beps(\dot\bv(s,\cdot))) \bI
\,,  \beps(\bv(s,\cdot)) \big) \dd s = 0.
\end{aligned}
\end{align}
Focusing in particular on the first and the third term on the left-hand side of \eqref{eq:47} we then have that
\begin{align}\label{eq:48}
\begin{aligned}
- \frac{1}{2} \tau^\alpha \int_0^{t_0} &\frac{\dd}{\dd s}\|\bu (s,\cdot)\|^2_{\LL^2_\varrho(\Omega)} \dd s  - (1-\tau^\alpha) \int_0^{t_0} ((-\dot{e}_{\alpha,1}\ast_s \varrho \dot{\bu})(s,\cdot), \bu(s,\cdot)) \dd s \\
& + \int_0^{t_0} \left( \frac{\dd}{\dd s}\|\beps(\bv(s,\cdot))\|^2_{\LL^2_\mu(\Omega)} + \frac{1}{2} \frac{\dd}{\dd s}\|\tr(\beps(\bv(s,\cdot)))\|^2_{\LL^2_\lambda(\Omega)}\right) \dd s = 0.
\end{aligned}
\end{align}
As
\begin{align*}
(-\dot{e}_{\alpha,1}\ast_s \dot{\bu})(s,\cdot) &= \frac{\dd}{\dd s} (-{e}_{\alpha,1}\ast_s \dot{\bu})(s,\cdot)
- (-{e}_{\alpha,1}(0)) \dot{\bu}(s,\cdot)
= \frac{\dd}{\dd s} (-{e}_{\alpha,1}\ast_s \dot{\bu})(s,\cdot)
+ \dot{\bu}(s,\cdot),
\end{align*}
inserting this into the second term on the left-hand side of \eqref{eq:48} yields
\begin{align*}
- \frac{1}{2}\tau^\alpha \int_0^{t_0} \frac{\dd}{\dd s}\|\bu (s,\cdot)\|^2_{\LL^2_\varrho(\Omega)} \dd s
&- (1-\tau^\alpha) \int_0^{t_0} \left(\frac{\dd}{\dd s} (-{e}_{\alpha,1}\ast_s \varrho\dot{\bu})(s,\cdot)+ \varrho \dot{\bu}(s,\cdot), \bu(s,\cdot)\right) \dd s \\
& + \int_0^{t_0} \left(\frac{\dd}{\dd s}\|\beps(\bv(s,\cdot))\|^2_{\LL^2_\mu(\Omega)} + \frac{1}{2} \frac{\dd}{\dd s}\|\tr(\beps(\bv(s,\cdot)))\|^2_{\LL^2_\lambda(\Omega)}\right) \dd s = 0.
\end{align*}
Hence, by performing partial integration in the second integral on the left-hand side, and because $\bv(t_0,\cdot)=0$, it follows that
\begin{align*}
&-\frac{1}{2}\tau^\alpha
\left(\|\bu(t_0,\cdot)\|^2_{\LL^2_\varrho(\Omega)}- \|\bu(0,\cdot)\|^2_{\LL^2_\varrho(\Omega)}\right)
- \frac{1}{2}(1-\tau^\alpha) \left(\|\bu(t_0,\cdot)\|^2_{\LL^2_\varrho(\Omega)} - \|\bu(0,\cdot)\|^2_{\LL^2_\varrho(\Omega)}\right)\\
& -  (1-\tau^\alpha)\int_0^{t_0} \left(({e}_{\alpha,1}\ast_s \sqrt{\varrho}\dot{\bu})(s,\cdot), \sqrt{\varrho}\dot{\bu}(s,\cdot)\right) \dd s
+ (1-\tau^\alpha) \left(({e}_{\alpha,1}\ast_s \varrho \dot{\bu})(s,\cdot), {\bu}(s,\cdot)\right)|_{s=0}^{s=t_0}\\
&- \|\beps(\bv(0,\cdot))\|^2_{\LL^2_\mu(\Omega)} - \frac{1}{2} \|\tr(\beps(\bv(0,\cdot)))\|^2_{\LL^2_\lambda(\Omega)} = 0.
\end{align*}
Again, because $\bu \in \CC([0,T];[\LL^2(\Omega)]^3)$ satisfies
$\bu(0,\bx) = \mathbf{0}$ for a.e. $\bx \in \Omega$, rearrangement yields
\begin{align}\label{eq:49}
\begin{aligned}
&\frac{1}{2} \|\bu(t_0,\cdot)\|^2_{\LL^2_\varrho(\Omega)} + (1-\tau^\alpha)\int_0^{t_0} \left(({e}_{\alpha,1}\ast_s \sqrt{\varrho} \dot{\bu})(s,\cdot), \sqrt{\varrho} \dot{\bu}(s,\cdot)\right) \dd s \\
& \qquad +  \|\beps(\bv(0,\cdot))\|^2_{\LL^2_\mu(\Omega)} + \frac{1}{2} \|\tr(\beps(\bv(0,\cdot)))\|^2_{\LL^2_\lambda(\Omega)}\\
& \qquad \qquad = (1-\tau^\alpha) \left(({e}_{\alpha,1}\ast_t \varrho \dot{\bu})(t_0,\cdot), {\bu}(t_0,\cdot)\right).
\end{aligned}
\end{align}
Thus, thanks to Lemma \ref{le:1} the second term on the left-hand side of \eqref{eq:49} can be bounded below, yielding
\begin{align}\label{eq:50}
\begin{aligned}
&\frac{1}{2} \|\bu(t_0,\cdot)\|^2_{\LL^2_\varrho(\Omega)} + \frac{1}{2} (1-\tau^\alpha)\left[  (e_{\alpha,1} \ast_t \|\dot\bu(\cdot)\|^2_{\LL^2_\varrho(\Omega)})(t_0)
+ \int_0^{t_0} e_{\alpha,1}(s)\|\dot\bu(s,\cdot)\|^2_{\LL^2_\varrho(\Omega)} \dd s \right]\\
& ~ ~ +  \|\beps(\bv(0,\cdot))\|^2_{\LL^2_\mu(\Omega)} + \frac{1}{2} \|\tr(\beps(\bv(0,\cdot)))\|^2_{\LL^2_\lambda(\Omega)}
\leq (1-\tau^\alpha) \left(({e}_{\alpha,1}\ast_t \varrho \dot{\bu})(t_0,\cdot), {\bu}(t_0,\cdot)\right).
\end{aligned}
\end{align}
Next, we will show that for any $t_0>0$ such that $t_0 \leq \min\left(T,1\right)$ the term on the right-hand side of \eqref{eq:50} can be completely absorbed into the left-hand side of the inequality.
Indeed, by Young's inequality, Minkowski's integral inequality, and the Cauchy--Schwarz inequality,
\begin{align*}
 (1-\tau^\alpha) \left(({e}_{\alpha,1} \right. &\ast_t \left.\varrho\dot{\bu})(t_0,\cdot), {\bu}(t_0,\cdot)\right) \\
 & \leq \frac{1}{2}  (1-\tau^\alpha)  \|\bu(t_0,\cdot)\|^2_{\LL^2_\varrho(\Omega)} + \frac{1}{2}  (1-\tau^\alpha) \left\| ({e}_{\alpha,1}\ast_t \dot{\bu})(t_0,\cdot)\right\|^2_{\LL^2_\varrho(\Omega)}\\
& \leq \frac{1}{2}  (1-\tau^\alpha)  \|\bu(t_0,\cdot)\|^2_{\LL^2_\varrho(\Omega)} + \frac{1}{2}   (1-\tau^\alpha) \left[({e}_{\alpha,1}\ast_t \|\dot{\bu}\|_{\LL^2_\varrho(\Omega)})(t_0)\right]^2\\
& \leq \frac{1}{2} (1-\tau^\alpha)  \|\bu(t_0,\cdot)\|^2_{\LL^2_\varrho(\Omega)}
+ \frac{1}{2} (1-\tau^\alpha) \left[\int_0^{t_0} {e}_{\alpha,1}(t_0 -s) \dd s\right] \left[({e}_{\alpha,1}\ast_t \|\dot{\bu}\|^2_{\LL^2_\varrho(\Omega)})(t_0)\right]\!.
\end{align*}
As $e_{\alpha,1}(0)=1$ and $t \in [0,\infty) \mapsto e_{\alpha,1}(t)$ is positive and monotonic decreasing, it follows that
\begin{align*}
(1-\tau^\alpha) \left(({e}_{\alpha,1}\ast_t \varrho \dot{\bu})(t_0,\cdot), {\bu}(t_0,\cdot)\right)
& \leq \frac{1}{2} (1-\tau^\alpha)  \|\bu(t_0,\cdot)\|^2_{\LL^2_\varrho(\Omega)}
+ \frac{1}{2} (1-\tau^\alpha) t_0 \left[({e}_{\alpha,1}\ast_t \|\dot{\bu}\|^2_{\LL^2_\varrho(\Omega)})(t_0)\right].
\end{align*}
Substituting this into the right-hand side of \eqref{eq:50} and, because $\tau \in (0,1]$ and $t_0\leq 1$, yields
\begin{align}\label{eq:51}
\begin{aligned}
\frac{1}{2} \tau^\alpha  \|\bu(t_0,\cdot)\|^2_{\LL^2_\varrho(\Omega)} + \frac{1}{2} (1-\tau^\alpha)\left[
\int_0^{t_0} e_{\alpha,1}(s)\|\dot\bu(s,\cdot)\|^2_{\LL^2_\varrho(\Omega)} \dd s \right]\\
 \quad +  \|\beps(\bv(0,\cdot))\|^2_{\LL^2_\mu(\Omega)} + \frac{1}{2} \|\tr(\beps(\bv(0,\cdot)))\|^2_{\LL^2_\lambda(\Omega)}
\leq 0.
\end{aligned}
\end{align}
Thus we deduce that $\bu(t,\cdot) = \mathbf{0}$ for all $t \in [0,t_0]$, for any $t_0>0$ such that $t_0 \leq \min(T,1)$.
If $t_0<T$, then having shown that $\bu(t,\cdot) = \mathbf{0}$ for all $t \in [0,t_0]$ we repeat the argument on successive time intervals $[kt_0, \min(T, (k+1)t_0)]$, $k=1,2,\dots, K$, with initial data $\bu(kt_0, \cdot)=\mathbf{0}$, $\bu_t(kt_0,\cdot) = \mathbf{0}$, where $K$ is the (unique) positive integer such that $Kt_0<T$ and $(K+1)t_0 \geq T$. Hence, $\bu(t,\cdot)=
\mathbf{0}$ for all $t \in [0,T]$. Thus we have shown the uniqueness of the weak solution.

\smallskip

STEP 5. \textit{Continuous dependence of the solution on the data.} As the problem under consideration is linear, the energy inequality \eqref{EE}
implies continuous dependence of weak solutions on the initial data and the load vector.

\smallskip

{\color{black}
STEP 6. \textit{Attainment of the initial condition for $\bsig$.}
By \eqref{eq:stress-defin} and noting that
\[ \mathcal{L}^{-1}\left(\frac{1 + p^\alpha}{1 + \tau^{\alpha} p^\alpha}\right) =  (\tau^{-\alpha}-1)\, \dot{e}_\alpha(t,\tau^{-\alpha}) + \tau^{-\alpha} \delta
\quad \mbox{ and } \quad \mathcal{L}^{-1}\left(\frac{p^{\alpha-1}}{1  + \tau^\alpha p^\alpha}\right)  = \tau^{-\alpha}\, e_\alpha(t,\tau^{-\alpha}), \]
we have that
\begin{align}\label{eq:sig-equa}
\begin{aligned}
\tau^\alpha \bsig(t,\cdot)& =  (1 -\tau^{\alpha}) \dot{e}_\alpha(t,\tau^{-\alpha}) \ast_t (2\mu \beps(\bu(t,\cdot)) + \lambda \tr(\beps(\bu(t,\cdot))) \bI)\\
&\quad + (2\mu \beps(\bu(t,\cdot)) + \lambda \tr(\beps(\bu(t,\cdot))) \bI)
+ e_\alpha(t,\tau^{-\alpha}) \,(\tau^\alpha\bs-
2\mu \beps(\bg) - \lambda \tr(\beps(\bg)) \bI)\\
& =: \bA_1 + \bA_2 + \bA_3.
\end{aligned}
\end{align}
We begin by showing that $\bsig$ belongs to $\CC_w([0,T];[\LL^2(\Omega)]^{3 \times 3})$.

As $e_\alpha(\cdot,\tau^{-\alpha}) \in \CC([0,\infty))$ and $\tau^\alpha\bs-2\mu \beps(\bg) - \lambda \tr(\beps(\bg)) \bI \in [\LL^2(\Omega)]^{3 \times 3}$ thanks to
\eqref{eq:19}, we have that $e_\alpha(\cdot,\tau^{-\alpha}) \,(\tau^\alpha\bs-
2\mu \beps(\bg) - \lambda \tr(\beps(\bg)) \bI)$ belongs to the function space $\CC([0,T];[\LL^2(\Omega)]^{3 \times 3})$, and therefore
also to  $\CC_w([0,T];[\LL^2(\Omega)]^{3 \times 3})$, implying that $\bA_3 \in \CC_w([0,T];[\LL^2(\Omega)]^{3 \times 3})$.

To show that $\bA_2 \in \CC_w([0,T];[\LL^2(\Omega)]^{3 \times 3})$, recall that $\bA_2 \in \LL^\infty(0,T;[\LL^2(\Omega)]^{3 \times 3})$, because $\bu$, as a weak solution, belongs to $\LL^\infty(0,T;[\HH^1_0(\Omega)]^{3})$. Together with the fact
that
\begin{equation}\label{eq:A2}
\dot{\bA}_2 \in \LL^\infty(0,T;[\HH^{-1}(\Omega)]^{3 \times 3}),
\end{equation}
which we shall now show, and the continuous and dense embedding of $[\LL^2(\Omega)]^{3 \times 3}$ into $[\HH^{-1}(\Omega)]^{3 \time 3}$, this will then yield that $\bA_2 \in \CC_w([0,T];[\LL^2(\Omega)]^{3 \times 3})$
(cf., again, eq. (8.49) in Lemma 8.1, Ch. 3 of \cite{LM}), as desired. To show that \eqref{eq:A2} holds, we appeal to the
following result from the theory of Sobolev spaces of Banach-space-valued functions (cf., for example, Theorem 1.4.40 on p.15 in \cite{CH}):

\smallskip

\textit{Suppose that} $X$ \textit{is a reflexive Banach space}, $I$ \textit{is a nonempty bounded open interval of} $\mathbb{R}$, \textit{and}  $u \in \LL^p(I;X)$ \textit{for some} $p \in [1,\infty]$. \textit{Then,} $u \in \WW^{1,p}(I;X)$ \textit{if, and only if, there exists a function}
$g \in  \LL^p(I;\mathbb{R})$ \textit{such that}
\[ \|u(t) - u(s)\|_X \leq  \left|\int_s^t g(\tau) \dd \tau\right|\]
\textit{for almost all $s, t \in I$, i.e., for all $s, t$ outside a common null set.}

\smallskip

We shall apply this result with $p=\infty$, $X = [\HH^{-1}(\Omega)]^{3 \times 3}$, and $g(\tau)=\|\dot{\bu}(\tau,\cdot)\|_{\LL^2(\Omega)}$. Clearly,
\begin{align*} \|\bA_2(t)  - \bA_2(s)\|_X &= \|2\mu \beps(\bu(t,\cdot)- \bu(s,\cdot)) + \lambda \tr(\beps(\bu(t,\cdot)-\bu(s,\cdot))) \bI
\|_X \\
&= \left\|2\mu \beps\left(\int_s^t\dot{\bu}(\tau,\cdot)\dd \tau \right) + \lambda \tr\left(\beps\left(\int_s^t\dot{\bu}(\tau,\cdot)\dd \tau \right)\right) \bI\right\|_X\\
& \leq 2\mu \left\|\beps\left(\int_s^t\dot{\bu}(\tau,\cdot)\dd \tau \right)\right\|_X + \lambda \left\| \tr\left(\beps\left(\int_s^t\dot{\bu}(\tau,\cdot)\dd \tau \right)\right) \bI\right\|_X\\
& \leq 2\mu \left\|\int_s^t\dot{\bu}(\tau,\cdot)\dd \tau\right\|_{\LL^2(\Omega)} + 3\lambda
\left\|\int_s^t\dot{\bu}(\tau,\cdot)\dd \tau\right\|_{\LL^2(\Omega)}\\
& \leq (2\mu + 3 \lambda) \left| \int_s^t \|\dot{\bu}(\tau,\cdot)\|_{\LL^2(\Omega)}\dd \tau \right| < \infty\qquad \forall\, s, t \in [0,T],
\end{align*}
because $\dot{\bu} \in \LL^\infty(0,T;[\LL^2(\Omega)]^3)$, where we have used the bound $\|\beps(\bw)\|_X \leq \|\bw\|_{\LL^2(\Omega)}$ with $X = [\HH^{-1}(\Omega)]^{3\times 3}$.
Therefore, $\bA_2 \in \WW^{1,\infty}(0,T;[\HH^{-1}(\Omega)]^{3\times 3})$, whereby also $\dot{\bA}_2 \in \LL^{\infty}(0,T;[\HH^{-1}(\Omega)]^{3\times 3})$.
Thus we have shown that $\bA_2 \in \CC_w([0,T];[\LL^3(\Omega)]^{3\times 3})$.

Concerning the term $\bA_1$, as $\bA_1 = (1 -\tau^{\alpha}) \dot{e}_\alpha(t,\tau^{-\alpha}) \ast_t \bA_2$, and
$\bA_2 \in \CC_w([0,T];[\LL^3(\Omega)]^{3\times 3})$, also $\bA_1 \in \CC_w([0,T];[\LL^3(\Omega)]^{3\times 3})$.

By summing $\bA_1$, $\bA_2$ and $\bA_3$ we thus deduce that $\bsig = \bA_1 + \bA_2 + \bA_3 \in \CC_w([0,T];[\LL^2(\Omega)]^{3 \times 3})$. It remains to prove the attainment of the initial condition
$\bsig(0,\cdot) = \bs(\cdot)$.

\bigskip

Thanks to Fubini's theorem and the continuity of the integral with respect to its (upper) limit,
\begin{align*}
&\lim_{t \rightarrow 0_+} (\dot{e}_\alpha(t,\tau^{-\alpha}) \ast_t (2\mu \beps(\bu(t,\cdot)) + \lambda \tr(\beps(\bu(t,\cdot))) \bI),\boldsymbol{W})\\
&\quad= \lim_{t \rightarrow 0_+} \int_0^t \dot{e}_\alpha(s,\tau^{-\alpha})\, ((2\mu \beps(\bu(t-s,\cdot)) + \lambda \tr(\beps(\bu(t-s,\cdot))) \bI),\boldsymbol{W}) \dd s = 0 \qquad \forall\, \boldsymbol{W} \in [\LL^2(\Omega)]^{3\times 3}.
\end{align*}
Hence, and noting that (recall that $\bA_2 \in \CC_w([0,T];[\LL^2(\Omega)]^{3 \times 3})$)
\[ \lim_{t \rightarrow 0_+} (2\mu \beps(\bu(t,\cdot)) + \lambda \tr(\beps(\bu(t,\cdot))) \bI, \boldsymbol{W}) =
  (2\mu \beps(\bg(\cdot)) + \lambda \tr(\beps(\bg(\cdot))) \bI, \boldsymbol{W})\qquad \forall \, \boldsymbol{W} \in [\LL^2(\Omega)]^{3 \times 3},\]
and because $e_\alpha(0,\tau^{-\alpha})=1$, we have from \eqref{eq:sig-equa} that
\[ \lim _{t \rightarrow 0_+} (\bsig(t,\cdot), \boldsymbol{W}(\cdot)) = (\bs(\cdot), \boldsymbol{W}(\cdot)) \qquad \forall \, \boldsymbol{W} \in [\LL^2(\Omega)]^{3 \times 3}.\]
Therefore, $\bsig(0,\cdot) = \bs(\cdot)$ as an equality in $\CC_w([0,T];[\LL^2(\Omega)]^{3 \times 3})$, as required.
\quad $\Box$
}

\bigskip
{\color{black}
The results of the paper can be straightforwardly extended to initial-boundary-value problems for the fractional Zener wave equation with mixed homogeneous Dirichlet/nonhomogeneous Neumann boundary
conditions, i.e., to problems where the domain boundary $\partial\Omega$ is the disjoint union of $\Gamma_{\rm D}$ and $\Gamma_{\rm N}$, with $\Gamma_{\rm D}$ having positive two-dimensional surface measure,
\begin{alignat*}{2}
\bu & = \boldsymbol{0} &&\qquad \mbox{on $\Gamma_{\rm D}$},\\
 [(2\mu \beps(\bu) + \lambda \tr(\beps(\bu)) \bI) \,+\, e_{\alpha,1}\, (\tau^\alpha\bs-
2\mu \beps(\bg) - \lambda \tr(\beps(\bg)) \bI)]\cdot \boldsymbol{n} & = \boldsymbol{s}&&\qquad \mbox{on $\Gamma_{\rm N}$},
\end{alignat*}
where $\boldsymbol{n}$ is the unit outward normal vector to $\partial\Omega$, and $\boldsymbol{s} \in \LL^\infty(0,T; [\LL^2(\Gamma_N)]^3)$ is given, at the expense of adding a term of the form
\[ \int_0^T \int_{\Gamma_{\rm N}} \boldsymbol{s}(t,\boldsymbol{\xi})\cdot \bv(t,\boldsymbol{\xi}) \dd \boldsymbol{\xi} \dd t \]
to the right-hand side of \eqref{eq:21}, replacing the function space $[\HH^1_0(\Omega)]^3$ throughout by the function space $[\HH^1_{\Gamma_{\rm D},0}(\Omega)]^3$ consisting of all functions in $[\HH^1(\Omega)]^3$ with zero trace on $\Gamma_{\rm D}$, and $[\HH^{-1}(\Omega)]^3$ signifying the dual space of $[\HH^1_{\Gamma_{\rm D},0}(\Omega)]^3$. In the special case when the initial stress $\bs$ is such that $\tau^\alpha\bs= 2\mu \beps(\bg) +\lambda \tr(\beps(\bg)) \bI$, the Neumann boundary condition on $\Gamma_{\rm N}$ and the source term $\bG$ in \eqref{eq:21}, defined by \eqref{eq:22}, are both simplified.
}

As a possible further, but now nontrivial, extension of the model \eqref{eq:4},  we note that Freed and Diethelm \cite{FD} have extended Fung's nonlinear constitutive law for soft biological tissues
into a constitutive law involving fractional time-derivatives in the
sense of Caputo, first in one space dimension and then in three space-dimensions.
The model is derived in a configuration that differs from the
current configuration by a rigid-body rotation; it being the polar configuration.
Freed and Diethelm introduce mappings for the fractional-order operators of integration
and differentiation between the polar and spatial configurations. They then use these
mappings in the construction of their proposed viscoelastic model. The mathematical analysis of
the associated set of partial differential equations, {\color{black} and the study of wave propagation governed by
the associated nonlinear system of nonlocal evolution equations} are beyond the scope of the present paper.

\smallskip

\textbf{Acknowledgements:} We are grateful to Professor Du\v{s}an Zorica (Mathematical Institute of the Serbian Academy of Sciences and Arts) for stimulating discussions, and to Sr{\dj}an Lazendi\'c (University of Ghent)
for drawing our attention to reference \cite{Siskova}.
Ljubica Oparnica is supported by the Serbian Ministry of Education, Science, and Technological Development, under the
grants 174005 and 174024, and by the FWO Odysseus project of Michael Ruzhansky.



\noindent
Ljubica Oparnica, \\
Faculty of Education, University of Novi Sad, Serbia

\smallskip

\noindent
\textit{\&}

\smallskip

\noindent
Department of Mathematics: Analysis, Logic and Discrete Mathematics,\\
University of Gent,
Krijgslaan 281, S8,
9000 Gent, Belgium\\
\texttt{\footnotesize Oparnica.Ljubica@UGent.be}

\bigskip

\noindent
Endre S\"uli,\\
Mathematical Institute,
University of Oxford,
Woodstock Road,
Oxford OX2 6GG, UK\\
\texttt{\footnotesize endre.suli@maths.ox.ac.uk}

\end{document}